\documentclass[a4paper,epic,eepic,leqno,11pt]{amsart}
\usepackage{amsmath,latexsym,amsfonts,amssymb,amscd,euscript,epic,eepic,times}

\numberwithin{equation}{section}

\newtheorem {theorem}[equation]         {Theorem}

\newtheorem {lemma}[equation]           {Lemma}
\newtheorem {proposition}[equation]     {Proposition}

\newtheorem {corollary}[equation]       {Corollary}

\theoremstyle{Definition}
\newtheorem {definition}[equation]      {Definition}
\newtheorem {remark}[equation]          {Remark}

\newtheorem {remarks}[equation]         {Remarks}

\makeatletter
\newcommand{\bfsubsection}
{\@startsection{subsection}{3}{0pt}{\baselineskip}
 {-\fontdimen2\font}{\normalfont\bfseries}}
\newcommand{\bfsubsubsection}
{\@startsection{subsubsection}{3}{0pt}{\baselineskip}
 {-\fontdimen2\font}{\normalfont\bfseries}}
\makeatother

\newfont{\bfc}{cmbsy10 scaled 1200}  
\newfont{\dr}{msbm10 scaled \magstep1}  
\newfont{\sdr}{msbm8}  
\newfont{\gl}{eufm10 scaled \magstep1}  

\newcommand{\Met}{\mbox{$Met$}}

\DeclareFontFamily{OT1}{rsfs}{}
\DeclareFontShape{OT1}{rsfs}{n}{it}{<->rsfs10}{}
\DeclareMathAlphabet{\curly}{OT1}{rsfs}{n}{it}

 \newcommand{\CC}{{\mathbb C}}

 \newcommand{\PP}{{\mathbb P}}
 
 \newcommand{\RR}{{\mathbb R}}

 \newcommand{\ZZ}{{\mathbb Z}}

 \newcommand{\PU}{{\textnormal{PU}}}

 \newcommand{\bmC}{\mbox{\boldmath$C$}}

 \newcommand{\bmL}{\mbox{\boldmath$L$}}
 \newcommand{\bmM}{\mbox{\boldmath$M$}}

 \newcommand{\bsigma}{{\boldsymbol{\sigma}}}

 \newcommand{\btau}{{\boldsymbol{\tau}}}

 \newcommand{\cA}{{\mathcal A}}

 \newcommand{\cE}{{\mathcal E}}
 \newcommand{\cF}{{\mathcal F}}

 \newcommand{\cK}{{\mathcal K}}
 
 \newcommand{\cM}{{\mathcal M}}
 
 \newcommand{\cO}{{\mathcal O}}

 \newcommand{\cR}{{\mathcal R}}

 \newcommand{\bcF}{{\boldsymbol{\mathcal{F}}}}

\newcommand{\AAA}{{\curly A}}

\newcommand{\FFF}{{\curly F}}
\newcommand{\GGG}{{\curly G}}

\newcommand{\NNN}{{\curly N}}

\newcommand{\RRR}{{\curly R}}

\newcommand{\euscA}{{\EuScript A}}

\newcommand{\euscM}{{\EuScript M}}

\newcommand{\euscS}{{\EuScript S}}

\newcommand{\YMH}{{Y\! M\! H}}
\newcommand{\qu}{/\kern-.7ex/}
\newcommand{\exh}{\to\kern-1.8ex\to}

 \newcommand{\ra}{\rightarrow}
 \newcommand{\lra}{\longrightarrow}
 \newcommand{\hra}{\hookrightarrow}

 \newcommand{\wt}{\widetilde} 

 \newcommand{\kahler}{K\"{a}hler}

\newcommand{\du}{\operatorname{d}}
\newcommand{\coproduct}{\operatorname{\amalg}}

\newcommand{\imag}{\mathop{{\fam0 \sqrt{-1}}}\nolimits}

\newcommand{\ad}{\operatorname{ad}}

\newcommand{\dbar}{\bar{\partial}}

\newcommand{\End}{\operatorname{End}}
\newcommand{\GL}{\operatorname{GL}}
  
\newcommand{\Hom}{\operatorname{Hom}}
  
\newcommand{\id}{\operatorname{id}}

\newcommand{\Lie}{\operatorname{Lie}}

\newcommand{\rk}{\operatorname{rk}}
\newcommand{\SL}{\operatorname{SL}}

\newcommand{\tr}{\operatorname{tr}}
\newcommand{\U}{\operatorname{U}}

\newcommand{\Vol}{\operatorname{Vol}}

 \newcommand{\JDG}{J. Diff. Geom.}

 \newcommand{\kah}{\kahler}

 \newcommand{\holfil}{{0\hra \cF_0\hra \cF_1\hra\cdots\hra  \cF_{m}=\cF}}

        {\begin{list}{}{%
                \settowidth{\labelwidth}{{\emph{#1:}}}%
                \setlength{\leftmargin}{\labelwidth+\labelsep}}}%
        {\end{list}}
\newenvironment{entry}
        {\begin{list}{}%
                {%
                  \setlength{\labelsep}{8mm}%
                  \setlength{\labelwidth}{40mm}%
                  \setlength{\leftmargin}{54mm}%
                  \setlength{\rightmargin}{5pt}%
                }%
        }%
        {\end{list}}
\newlength{\Mylen}
\newcommand{\Lentrylabel}[1]{%
        \settowidth{\Mylen}{\emph{#1}}%
        \ifthenelse{\lengthtest{\Mylen > \labelwidth}}%
                {\parbox[b]{\labelwidth}
                        {\makebox[0pt][l]{\emph{#1}}\\}}%
                {\emph{#1}}
        \hfil\relax}

\newcommand{\Mentrylabel}[1]%
        {\raisebox{0pt}[1ex][0pt]{\makebox[\labelwidth][l]%
                {\parbox[t]{\labelwidth}{\hspace{0pt}\emph{#1:}}}}}
        {\begin{entry}}%
        {\end{entry}}

  {\begin{list}{}{%
   \settowidth{\labelwidth}{{#1}}%
   \setlength{\leftmargin}{\labelwidth+\labelsep}}}
  {\end{list}}

\title[Hitchin--Kobayashi correspondence, quivers, and vortices]
{Hitchin--Kobayashi correspondence, quivers, and vortices}
\author{Luis \'Alvarez--C\'onsul}  
\address{Department of Mathematics, University of Illinois at
Urbana-Champaign, Urbana, IL 61801, USA}
\curraddr{Mathematical Sciences, University of Bath, Bath, BA2 7AY, UK}
\email{L.Alvarez-Consul@maths.bath.ac.uk}
\author{Oscar Garc\'{\i}a--Prada}
\address{Departamento de Matem{\'a}ticas \\
Universidad Aut{\'o}noma de Madrid \\
28049 Madrid, SPAIN}
\curraddr{Instituto de Matem\'aticas y F\'{\i}sica Fundamental, CSIC, 
Serrano 113 bis, 28006 Madrid, Spain}
\email{oscar.garcia-prada@uam.es}

\begin{document}
\maketitle
\begin{abstract}
A twisted quiver bundle is a set of holomorphic vector 
bundles over a complex manifold, labelled by the vertices of a quiver,
linked by a set of morphisms twisted by a fixed collection of
holomorphic vector bundles, labelled by the arrows. When the manifold
is \kah , quiver bundles admit natural gauge-theoretic 
equations, which unify many 
known equations for bundles with extra structure. In this 
paper we prove a Hitchin--Kobayashi correspondence for 
twisted quiver bundles over a compact \kah\ manifold, relating the 
existence of solutions to the gauge equations to a stability criterion, 
and consider its application to a number of situations related to Higgs 
bundles and dimensional reductions of the Hermitian--Einstein 
equations.
\end{abstract}

\section*{Introduction}

A quiver $Q$ consists of a set $Q_0$ of vertices $v,v',\ldots$, and a
set $Q_1$ of arrows $a:v\to v'$ connecting the vertices.
Given a quiver  and a compact K\"ahler manifold, 
a {\em quiver bundle} is defined by assigning a holomorphic vector 
bundle $\cE_v$ to a finite number of vertices and a homomorphism 
$\phi_a:\cE_v\to\cE_{v'}$ to a finite number of  
arrows. A {\em quiver sheaf} is defined by replacing the term 
`holomorphic vector bundle' by `coherent sheaf' in this definition.
If we fix a collection of holomorphic vector bundles $M_a$ 
parametrized by the set of arrows, 
and the morphisms are $\phi_a:\cE_v\otimes M_a\to\cE_{v'}$, 
twisted by the corresponding bundles, we have a 
{\em twisted quiver bundle} or a {\em twisted quiver sheaf}. 
In this paper we define natural gauge-theoretic equations, that we call 
{\em quiver vortex equations}, for a collection of hermitian 
metrics on the bundles associated to the vertices of a twisted quiver 
bundle (for this, we need to fix hermitian metrics on the twisting 
vector bundles). To solve these equations, we introduce a stability criterion 
for twisted quiver sheaves, and prove a Hitchin--Kobayashi correspondence, 
relating the existence of (unique) hermitian metrics satisfying the quiver 
vortex equations to the stability of the quiver bundle. 
The equations and the stability criterion depend 
on some real numbers, the {\em stability parameters} (cf. Remarks \ref{remark:stability} 
for the exact number of parameters). 
It is relevant to point out that our results cannot be derived from the 
general Hitchin--Kobayashi correspondence scheme developed
by Banfield \cite{Ba} and further generalized by Mundet \cite{M}. 
This is due not only to the presence of twisting vector bundles, but also
to the deformation of the Hermitian--Einstein terms in the equations. 
This deformation  is naturally explained by the symplectic interpretation  
of the equations, and  accounts for extra parameters in 
the stability condition for the twisted quiver bundle. 

This correspondence provides a unifying framework to study a number of 
problems that have been considered 
previously. 
The simplest situation occurs when the quiver has a single vertex and
no arrows, in which case a quiver bundle is just a holomorphic bundle
$\cE$, and the gauge equation is the {\em Hermitian--Einstein
  equation}. A theorem of Donaldson, Uhlenbeck and Yau
\cite{D1,D2,UY}, establishes that a (unique) solution to the
Hermitian--Einstein equation exists if and only if $\cE$ is
polystable. The bundle $\cE$ is called {\em stable}  
(in the sense of Mumford--Takemoto) if 
$\mu(\cF)<\mu(\cE)$ for each proper coherent subsheaf $\cF\subset
\cE$, where the slope $\mu(\cF)$ 
is the degree divided by the rank; a finite direct sum of 
stable bundles with the same slope is called {\em polystable}. 
A correspondence of this type is usually  known as a {\em 
Hitchin--Kobayashi correspondence}.
A Hitchin--Kobayashi correspondence, where some extra 
structure is added to the bundle $\cE$, appears in the theory of 
{\em Higgs bundles}, consisting of 
pairs $(\cE,\Phi)$ formed by a holomorphic vector bundle $\cE$ 
and a morphism $\Phi:\cE\to\cE\otimes\Omega$, where $\Omega$ is the 
sheaf of holomorphic differentials (sometimes the condition 
$\Phi\wedge\Phi=0$ is added as part of the definition). Higgs bundles
were first studied by Hitchin \cite{H} (when $X$ is a compact Riemann
surface), and Simpson \cite{S} (when $X$ is higher dimensional), who
introduced a natural gauge equation for them, and proved a
Hitchin--Kobayashi correspondence. Higgs bundles are twisted quiver
bundles, for a quiver formed by one vertex and one arrow whose head
and tail coincide, and the twisting bundle is the holomorphic tangent
bundle (i.e. the dual to $\Omega$). 
Another class of quiver bundles are {\em holomorphic triples} 
$(\cE_1,\cE_2,\Phi)$, consisting of two holomorphic 
bundles $\cE_1$ and $\cE_2$, and a morphism $\Phi:\cE_2\to\cE_1$. 
The quiver has two vertices, say $1$ and $2$, and one arrow $a:2\to 1$ 
(the twisting sheaf is $\cO_X$). The corresponding equations 
are called the coupled vortex equations \cite{G2,BG}.
When $\cE_2=\cO_X$, holomorphic triples are {\em holomorphic pairs}
$(\cE,\Phi)$, where $\cE$ is a bundle and $\Phi\in H^0(X,\cE)$
(cf. \cite{B}). 

There are other examples of  quiver vortex equations that
come out naturally from the study of the
moduli of solutions to the Higgs bundle equation. 
Combining a theorem of Donaldson and Corlette \cite{D3,C} with 
the Hitchin--Kobayashi correspondence for Higgs bundles \cite{H,S}, 
one has that the set of isomorphism classes of semisimple 
complex representations of the fundamental group of $X$ in 
$\GL(r,\CC)$ is in bijection with the moduli space 
of polystable Higgs bundles with vanishing Chern classes.
When $X$ is a compact Riemann surface, this generalizes 
a theorem of Narasimhan and Seshadri \cite{NS}, which provides an 
interpretation of the unitary representations of the fundamental 
group as degree zero polystable vector bundles, up to isomorphism. 
Now, if $X$ is a compact Riemann surface of genus $g\geq 2$, the
Morse methods introduced by Hitchin \cite{H} reduce the study of the topology
of the moduli $\cM$ of Higgs bundles to the study of the topology
of the moduli of complex variations of the Hodge structure --- 
the critical points of the Morse function in this case.
These are twisted quiver bundles, called {\em twisted holomorphic
chains}, for a quiver whose vertex set is the set $\ZZ$ of integer
numbers, and whose arrows are $a_i:i\to i+1$, for each $i\in\ZZ$; the
twisting bundle associated to each arrow is the holomorphic tangent
bundle. The twisted holomorphic chains that appear in these critical 
submanifolds are polystable for particular values of the stability
parameters. Using Morse theory, Hitchin \cite{H} computed the
Poincar\'e polynomial of $\cM$ for the rank 2 case. Gothen \cite{Go}
obtained similar results for rank 3: the critical submanifolds are
moduli spaces of stable twisted holomorphic chains formed by a line
bundle and a rank 2 bundle (i.e. twisted holomorphic triples), and by
three line bundles. To use these methods for higher rank, one needs
to study moduli spaces of other twisted holomorphic chains. 
A possible strategy is to proceed as in \cite{Th}, 
studying the moduli of twisted holomorphic chains in the whole 
parameter space. Another interesting type of quiver bundles arise in 
the study of semisimple representations of the fundamental group of $X$
in $\U(p,q)$, the unitary group for a hermitian 
inner product of indefinite signature.
Here, the quiver has two vertices, say
$1$ and $2$, and two arrows, $a:1\to 2$ and $b:2\to 1$, and the twisting
bundle associated to each arrow is the holomorphic tangent bundle. 
These are  studied in \cite{BGG1,BGG2}. 

Another context in which quiver bundles 
appear naturally is in the study of dimensional reductions
 of the Hermitian--Einstein 
equation over the product of a \kah\ manifold $X$ and a 
flag manifold. In this case, the parabolic subgroup defining the flag
manifold  entirely determines the
structure of the quiver \cite{AG1,AG2}. The dimensional reduction 
for this kind of manifolds has provided insight in the general theory of quiver 
bundles, and was actually the first method used to prove a Hitchin--Kobayashi
correspondence for holomorphic triples 
\cite{G2,BG}, holomorphic chains \cite{AG1}, and quiver bundles for more 
general quivers with relations \cite{AG2}. In these 
examples, the quiver bundles are not twisted, however,  there are other examples 
for which a generalization of the method of dimensional reduction has 
produced {\em twisted holomorphic triples} \cite{BGK1,BGK2}. 

An important feature of the stability of quiver sheaves is that 
it generally depends on  several real parameters.
When $X$ is an algebraic variety, the ranks and degrees 
appearing in the numerical condition defining the stability criterion 
are integral, and the parameter space is partitioned 
into chambers. Strictly semistable quiver sheaves can occur 
when the parameters are on a wall separating the chambers, and 
the stability condition 
only depends on the chamber in which the parameters are. 
In the case of holomorphic triples \cite{BG}, there 
is a chamber (actually an interval in $\RR$) where the stability of the
triple is related to the stability of the bundles. This can be used to
obtain existence theorems for stable triples when the  
parameters are in this chamber, while the methods of \cite{Th} can 
be used to prove existence results for other chambers 
(see \cite{BGG2} for recent work in the case of triples). 
The geography of the resulting convex polytope for other 
quivers is an interesting issue to which we wish to return in a 
future paper. To approach this problem, one should study the 
homological algebra of quiver bundles. This has been developed by 
Gothen and King in a paper \cite{GK} that appeared after we submitted
this paper.

When the manifold $X$ is a point, a quiver bundle is just a quiver
module (over $\CC$; cf. e.g. \cite{ARS}). For
arbitrary $X$, a quiver bundle can be regarded as a family of quiver modules
(the fibres of the quiver bundle), parametrized by $X$. One can thus 
transfer to our setting many constructions of the theory of quiver
modules. In the last part of the paper we introduce a more algebraic point of view
by considering the path algebra bundle of the twisted quiver and looking 
at twisted quiver bundles as locally free modules over this
bundle of algebras. This point of view is inspired by a similar
construction for quiver modules \cite{ARS}, and suggests a 
generalization to other algebras that appear naturally in other
problems. This is something to which we plan to come back in the future. 

The Hitchin--Kobayashi correspondence for quiver bundles 
combines in one theory two different versions, in some sense, 
of the theorem 
of Kempf and Ness \cite{KN} 
identifying the symplectic quotient of a projective variety
by a compact Lie group action, with the geometric invariant 
theory quotient. The first one is
the classical Hitchin--Kobayashi correspondence for vector bundles, 
and the second one occurs when the
manifold $X$ is a point, in which case the equations and the 
stability condition reduce to the moment 
map equations and the stability condition for quiver modules 
introduced by King \cite{K}. 
As we prove in Theorem \ref{thm:Bogomolov}, 
there is in fact a very tight relation between the  quiver vortex 
equations and  the moment map 
equations for quiver modules: when the twisting sheaves 
are $\cO_X$ and the bundles have vanishing Chern classes, the
existence of solutions to the quiver vortex equations is equivalent to
the existence of flat metrics on the bundles which fibrewise satisfy
the moment map equations for quiver modules. 

\section{Twisted quiver bundles}

In this section we define the basic objects that we shall study: {\em
twisted quiver bundles} and {\em twisted quiver sheaves}. They
are representations of quivers in the categories of holomorphic vector
bundles and coherent sheaves, respectively, twisted by some fixed
holomorphic vector bundles, as explained in \S \ref{sub:quiverbundles}. 
Thus, many results about quiver modules, i.e. quiver
representations in the category of vector spaces, can be tranferred to
our setting. A good reference for quivers and their linear
representations is \cite{ARS}. 

\subsection{Quivers}
\label{sub:quivers}

A quiver, or directed graph, is a pair of sets $Q=(Q_0,Q_1)$
together with two maps $h,t:Q_1\ra Q_0$. The elements of $Q_0$
(resp. $Q_1$) are called the vertices (resp. arrows) of the
quiver. For each arrow $a\in Q_1$, the vertex $ta$ (resp. $ha$) is
called the tail (resp. head) of the arrow $a$. The arrow $a$ is 
sometimes represented by $a:v\to v'$ when $v=ta$ and $v'=ha$.

\subsection{Twisted quiver sheaves and bundles}
\label{sub:quiverbundles}

Throughout this paper, $X$ is a connected compact \kah\ manifold, $Q$
is a quiver, and $M$ is a collection of finite rank
locally free sheaves $M_a$ on $X$, for each arrow $a\in Q_1$. By a 
sheaf on $X$, we shall will mean an analytic sheaf of $\cO_X$-modules.
Our basic objects are given by the following: 

\begin{definition}
\label{def:Q-sheaf}
An {\em $M$-twisted $Q$-sheaf} on $X$ is a pair $\cR=(\cE ,\phi)$,
where $\cE$ is a collection of coherent sheaves $\cE_v$ on $X$, for
each $v\in Q_0$, and $\phi$ is a collection of morphisms 
$\phi_a:\cE_{ta}\otimes M_a\to\cE_{ha}$, for each $a\in Q_1$, such 
that $\cE_v=0$ for all but finitely many $v\in Q_0$, and 
$\phi_a=0$ for all but finitely many $a\in Q_1$.
\end{definition}

\begin{remark}
\label{remark:finite-quiver}
{\rm 
Given a quiver $Q=(Q_0,Q_1)$, as defined 
in \S \ref{sub:quivers}, the sets $Q_0$ and $Q_1$ can be infinite, 
but for each $M$-twisted $Q$-sheaf $\cR=(\cE,\phi)$, the subset 
$Q_0'\subset Q_0$ of vertices $v$ such that $\cE_v\neq 0$, and the 
subset $Q_1'\subset Q_1$ of arrows $a$ such that $\phi_a\neq 0$, are
both finite. Thus, to any $M$-twisted $Q$-sheaf $\cR=(\cE,\phi)$, 
we can associate the subquiver $Q'=(Q_0',Q_1')$ of $Q$, and $\cR$ 
can be seen as an $M'$-twisted $Q'$-sheaf, where $Q_0'$, $Q_1'$ 
are finite sets, and $M'\subset M$ is the collection 
of sheaves $M_a$ with $a\in Q_1'$.
}\end{remark}

As usual, we identify a holomorphic vector bundle $\cE$, with the locally
free sheaf of sections of $\cE$. Accordingly, a {\em holomorphic
$M$-twisted $Q$-bundle} is an $M$-twisted $Q$-sheaf $\cR=(\cE,\phi)$ such
that the sheaf $\cE_v$ is a holomorphic vector bundle, for each $v\in
Q_0$.  
For the sake of brevity, in the following the terms `$Q$-sheaf' or
`$Q$-bundle' are to be understood as `$M$-twisted $Q$-sheaf' or 
`$M$-twisted $Q$-bundle', respectively, often suppressing the adjective
`$M$-twisted'.

A {\em morphism} $f:\cR\ra\cR'$ between two
$Q$-sheaves $\cR=(\cE,\phi)$, $\cR'=(\cE',\phi')$, is given by a
collection of morphisms $f_v:\cE_v\to\cE'_v$, for each $v\in Q_0$, such
that $\phi'_a\circ (f_v\otimes \id_{M_a})=f_{v'}\circ\phi_a$, for each arrow $a:v\to v'$
in $Q$. If $f:\cR\ra\cR'$ and $g:\cR'\ra\cR''$ are two
morphisms between representations $\cR=(\cE,\phi)$, $\cR'=(\cE',\phi')$, $\cR''
=(\cE'',\phi'')$, then the composition $g\circ f$ is defined as
the collection of composed morphisms $g_v\circ f_v:\cE_v\to\cE''_v$, for each $v\in
Q_0$. We have thus defined the category of $M$-twisted $Q$-sheaves on $X$,
which is abelian. Important concepts in relation to
stability and semistability (defined in \S \ref{sub:Q-stability}) are the
notions of $Q$-subsheaves and quotient $Q$-sheaves, as well as
indecomposable and simple $Q$-sheaves. They are defined as for any
abelian category. In particular, an $M$-twisted $Q$-subsheaf of
$\cR=(\cE,\phi)$ is another $M$-twisted $Q$-sheaf
$\cR'=(\cE',\phi')$ such that $\cE'_v\subset\cE_v$,
for each $v\in Q_0$, $\phi_a(\cE_{ta}'\otimes M_a)\subset\cE_{ha}'$,
for each $a\in Q_1$, and $\phi'_a:\cE'_{ta}\otimes M_a\to\cE'_{ha}$ is
the restriction of $\phi_a$ to $\cE_{ta}'\otimes M_a$, for each $a\in
Q_0$. 

\section{Gauge equations and stability}

\subsection{Gauge equations}
\label{sub:vortex-equation}

Throughout this paper, given a smooth bundle $E$ on $X$,
$\Omega^{k}(E)$ (resp. $\Omega^{i,j}(E)$) is the space of smooth
$E$-valued complex $k$-forms (resp. $(i,j)$-forms) on $X$, $\omega$
is a fixed \kah\ form on $X$, and $\Lambda:\Omega^{i,j}(E)\to
\Omega^{i-1,j-1}(E)$ is contraction with $\omega$
(we use the same notation as e.g. in \cite{D1}). 
The gauge equations will also depend on a fixed
collection $q$ of hermitian metrics $q_a$ on $M_a$, for each $a\in
Q_1$, which we fix once and for all. Let $\cR=(\cE,\phi)$ 
be a holomorphic $M$-twisted $Q$-bundle on $X$. A {\em hermitian metric} on
$\cR$ is a collection $H$ of hermitian metrics $H_v$ on $\cE_v$,
for each $v\in Q_0$ with $\cE_v\neq 0$. 
To define the gauge equations on $\cR$, we note that
$\phi_a:\cE_{ta}\otimes M_a\to\cE_{ha}$ has a smooth adjoint morphism $\phi_a^{*H_a}:\cE_{ha}\to\cE_{ta}\otimes M_a$ with respect to the hermitian metrics $H_{ta}\otimes
q_a$ on $\cE_{ta}\otimes M_a$, and $H_{ha}$ on 
$\cE_{ha}$, for each $a\in Q_0$, so it makes sense to consider the composition 
$\phi_a\circ\phi^{*H_a}_a:\cE_{ha}\to \cE_{ta}\otimes M_a\to\cE_{ha}$. 
Moreover, $\phi_a$ and $\phi^{*H_a}_a$ can be seen as morphisms $\phi_a:
\cE_{ta}\to\cE_{ha} \otimes M^*_a$ and $\phi^{*H_a}_a:\cE_{ha}\otimes M^*_a\to 
\cE_{ta}$, so $\phi^{*H_a}_a\circ\phi_a:\cE_{ta}\to
\cE_{ta}$ makes sense too.

\begin{definition}\label{def:vortex-equation}
Let $\sigma$ and $\tau$ be collections of real numbers $\sigma_v,
\tau_v$, with $\sigma_v$ positive, for each $v\in Q_0$. 
A hermitian metric $H$ satisfies the $M$-twisted quiver
$(\sigma,\tau)$-vortex equations if 
\begin{equation}\label{eq:vortex-equation}
\sigma_v \imag \Lambda F_{H_v}
+\sum_{a\in h^{-1}(v)}\phi_a\circ\phi_a^{*H_a}-\sum_{a\in t^{-1}(v)}\phi_a^{*H_a}\circ\phi_a 
=\tau_v\id_{\cE_v}, \quad\quad 
\end{equation}
for each $v\in Q_0$ such that $\cE_v\neq 0$, where $F_{H_v}$ is the
curvature of the Chern connection $A_{H_v}$ associated to the metric
$H_v$ on the holomorphic vector bundle $\cE_v$, for each $v\in Q_0$ with
$\cE_v\neq 0$. 
\end{definition}

\subsection{Moment map interpretation}
\label{sub:momentmap}

The twisted quiver vortex equations appear as a symplectic 
reduction condition, as we explain now. Let $E$ be a collection of
smooth vector bundles $E_v$, for each $v\in Q_0$, with $E_v=0$ for all
but finitely many $v\in Q_0$. By removing the vertices $v\in Q_0$ with
$E_v=0$ and all but finitely many arrows $a\in Q_1$, we
obtain a finite subquiver, which we still call $Q=(Q_0,Q_1)$, such
that $E_v\neq 0$ for each $v\in Q_0$ 
(see Remark \ref{remark:finite-quiver}). Let $H_v$ be a hermitian metric
on $E_v$, for each $v\in Q_0$. Let $\AAA_v$ and $\GGG_v$ be the 
corresponding spaces of unitary connections and their unitary gauge
groups, and let $\AAA^{1,1}_v\subset\AAA_v$ be the space of unitary
connections $A_v$ with $(\dbar_{A_v})^2=0$, for each $v\in Q_0$. The
group 
$$
\GGG=\prod_{v\in Q_0}\GGG_v
$$
acts on the space $\AAA$ of unitary connections, and on the
representation space $\Omega^0$, defined by
\begin{equation}
  \label{eq:AAA-RRR}
  \AAA=\prod_{v\in Q_0}\AAA_v ,\quad \Omega^0=\Omega^0(\RRR(Q,E)),
{\rm ~~with~}
\RRR(Q,E)=\bigoplus_{a\in Q_1} \Hom(E_{ta}\otimes M_a,E_{ha}), 
\end{equation}
where $\Hom(E_{ta}\otimes M_a,E_{ha})$ is the vector bundle of
homomorphisms $E_{ta}\otimes M_a\to E_{ha}$. 
An element $g\in\GGG$ is a collection of group elements 
$g_v\in\GGG_v$, for each $v\in Q_0$, and an element
$A\in\AAA$ (resp. $\phi\in\Omega^0$) is a 
collection of unitary connections $A_v\in\AAA_v$
(resp. smooth morphisms $\phi_a:E_{ta}\otimes M_a\to E_{ha}$), for each
$v\in Q_0$ (resp. $a\in Q_1$). 
The $\GGG$-actions on $\AAA$ and $\Omega^0$ are
$
\GGG\times\AAA\to\AAA, (g,A)\mapsto A'=g\cdot A , 
{\rm ~with~} \du_{A'_v}=g_v\circ \du_{A_v} \circ g_v^{-1}, {\rm
~for~each~} v\in Q_0;
$
$\GGG\times\Omega^0\to\Omega^0, (g,\phi)\mapsto \phi'=g\cdot\phi , 
{\rm ~with~} \phi'_a=g_{ha} \circ\phi_a\circ (g_{ta}^{-1}\otimes\id_{M_a}), 
{\rm ~for~each~} a\in Q_1,$
respectively. 
The induced $\GGG$-action on the product $\AAA\times
\Omega^0$ leaves invariant the subset $\NNN$ of pairs
$(A,\phi)$ such that $A_{v}\in\AAA^{1,1}_v$,
for each $v\in Q_0$, and $\phi_a:E_{ta}\otimes M_a\to E_{ha}$ is
holomorphic with respect to $\dbar_{A_{ta}}$ and $\dbar_{A_{ha}}$, for
each $a\in Q_0$.
Let $\omega_v$ be the
$\GGG_v$-invariant symplectic form on $\AAA_v$, for each
$v\in Q_0$, as given in \cite{AB} for a compact Riemann surface, or
e.g. in \cite[Proposition 6.5.8]{DK} for any compact \kah\ manifold,
that is, 
$$
\omega_v(\xi_v,\eta_v)=\int_X\Lambda\tr(\xi_v\wedge\eta_v), \quad
{\rm for~} \xi_v,\eta_v\in\Omega^1(\ad (E_v)), 
$$
where $\ad (E_v)$ is the vector bundle of $H_v$-antiselfadjoint
endomorphisms of $E_v$. The corresponding moment map
$\mu_v:\AAA_v\to(\Lie\GGG_v)^*$ is given by $\mu_v(A_v)=\Lambda F_{A_v}$
(we use implicitly the inclusion of $\Lie\GGG_v$ in its dual space
by means of the metric $H_v$ on $E_v$). The symplectic form $\omega_{\RRR}$
on $\Omega^0$ associated to the $L^2$-metric induced by the 
hermitian metrics on the spaces $\Omega^0(\Hom(E_{ta}\otimes M_a,E_{ha}))$ 
is $\GGG$-invariant, and has associated moment map
$\mu_{\RRR}:\Omega^0\to(\Lie\GGG)^{*}$ given by 
$\mu_{\RRR}=\sum_{v\in Q_0}\mu_{\RRR,v}$, 
with $\mu_{\RRR, v}:\Omega^0\to\Lie\GGG_v\subset\Lie\GGG\subset
\Lie\GGG)^{*}$ given by  
\begin{equation}
  \label{eq:mmap-RRR-bis}
  \imag\mu_{\RRR , v}(\phi)=
  \sum_{a\in h^{-1}(v)}\phi_a\circ\phi_a^{*H_a}-\sum_{a\in
    t^{-1}(v)}\phi_a^{*H_a}\circ\phi_a, \quad {\rm for~} \phi\in\Omega^0,
\end{equation}
(this follows as in \cite[\S 6]{K}, which considers the action of
a unitary group on a representation space of quiver modules).
Given a collection $\sigma$ of real numbers $\sigma_v>0$, for
each $v\in Q_0$, $\sum_{v\in Q_0}\sigma_v\omega_v+\omega_{\RRR}$ is
obviously a $\GGG$-invariant symplectic form on $\AAA\times\Omega^0$.
A moment map for this symplectic form is
$\mu_{\sigma}=\sum_{v\in Q_0}\sigma_v\mu_v+\mu_{\RRR}$, 
where we are omitting pull-backs to $\AAA\times\Omega^0$ in the
notation. Any collection $\tau$ of real numbers $\tau_v$, for each
$v\in Q_0$ defines an element $\imag\tau\cdot\id=\imag \sum_{v\in 
Q_0} \tau_v\id_{\cE_v}$ in the center of $\Lie\GGG$. 
The points of the symplectic reduction $\mu_{\sigma}^{-1}
(-\imag\cdot\tau)/\GGG$ are precisely the orbits of 
pairs $(A,\phi)$ such that the hermitian metric $H$ satisfies the
$M$-twisted $(\sigma,\tau)$-vortex quiver equations on the corresponding
holomorphic quiver bundle $\cR=(\cE,\phi)$. 
Thus, Definition \ref{def:vortex-equation} 
picks up the points of $\mu_{\sigma }^{-1}(-\imag\tau)$ in the \kah\
submanifold (outside its singularities) $\NNN$. For convenience in the
Hitchin--Kobayashi correspondence, it is formulated in terms of
hermitian metrics.

\subsection{Stability}\label{sub:Q-stability}

To define stability, we need some preliminaries and notation. Let $n$
be the complex dimension of $X$. Given a torsion-free coherent sheaf
$\cE$ on $X$, the double dual sheaf 
$\det(\cE)^{**}$ is a holomorphic line bundle, and we define the first
Chern class $c_1(\cE)$ of $\cE$ as the first Chern class of
$\det(\cE)^{**}$. The degree of $\cE$ is the real number 
$$
\deg(\cE)=\frac{2\pi}{\Vol(X)}\frac{1}{(n-1)!}
\left< c_1(\cE)\smile\left[\omega^{n-1}\right],[X]\right> , 
$$
where $\Vol(X)$ is the volume of $X$, $[\omega^{n-1}]$ is the
cohomology class of $\omega^{n-1}$, and 
$[X]$ is the fundamental class of $X$. Note that the degree
depends on the cohomology class of $\omega$. Given a holomorphic
vector  bundle $\cE$ on $X$, by Chern-Weil theory, its degree equals 
$$
\deg(\cE)=\frac{1}{\Vol(X)}\int_X\tr(\imag\Lambda F_{H}),
$$
where $F_H$ is the curvature of the Chern connection associated
to a hermitian metric $H$ on $\cE$. 

Let $Q$ be a quiver, and $\sigma$, $\tau$ be
collections of real numbers $\sigma_v, \tau_v$, with $\sigma_v>0$, 
for each $v\in Q_0$; $\sigma$ and $\tau$ are called the {\em stability
parameters}. Let $\cR=(\cE,\phi)$ be a $Q$-sheaf on $X$. 

\begin{definition}\label{def:stab-Q-sheaf}
The $(\sigma,\tau)$-{\em degree} and $(\sigma,\tau)$-{\em slope}
of $\cR$ are
$$
\deg_{\sigma ,\tau }(\cR)=\sum_{v\in Q_0}
\left(\sigma_v\deg(\cE_v)-\tau_v\rk(\cE_v)\right) , 
\quad
\mu_{\sigma ,\tau }(\cR)=\frac{\deg_{\sigma ,\tau }(\cR)}{\sum_{v\in Q_0}
\sigma_v\rk(\cE_v)}, 
$$
respectively. The $Q$-sheaf $\cR$ is called $(\sigma ,\tau)$-{\em
stable} (resp. $(\sigma ,\tau)$-{\em semistable}) if for all proper
$Q$-subsheaves $\cR'$ of $\cR$, $\mu_{\sigma ,\tau}(\cR')<\mu_{\sigma
,\tau }(\cR)$ (resp. $\mu_{\sigma ,\tau}(\cR')\leq\mu_{\sigma ,\tau }
(\cR)$). A $(\sigma,\tau)$-{\em polystable} $Q$-sheaf is a direct sum
of $(\sigma,\tau)$-stable $Q$-sheaves, all of them with the same
$(\sigma,\tau)$-slope. 
\end{definition}

As for coherent sheaves, one can prove that any $(\sigma ,\tau)$-stable
$Q$-sheaf is simple, i.e. its only endomorphisms are the multiples of
the identity.

\begin{remarks}\label{remark:stability}{\rm 
\begin{enumerate}
\item[(i)] 
If a holomorphic $Q$-bundle $\cR$ 
admits a hermitian metric satisfying the $(\sigma,\tau)$-vortex
equations, then taking traces in \eqref{eq:vortex-equation},
summing for $v\in Q_0$, and integrating over $X$, we see that the
parameters $\sigma,\tau$ are constrained by $\deg_{\sigma ,\tau
}(\cR)=0$. 
\item[(ii)]
If we transform the parameters $\sigma,\tau$, multiplying by a global
constant $c>0$, obtaining $\sigma'=c\sigma$, $\tau'=c\tau$, then
$\mu_{\sigma', \tau'}(\cR) =\mu_{\sigma ,\tau }(\cR)$. Furthermore, if
we transform the parameters $\tau$ by $\tau'_v=\tau_v+d\sigma_v$ for some
$d\in\RR$, and let $\sigma'=\sigma$, then $\mu_{\sigma',\tau'}(\cR) 
=\mu_{\sigma ,\tau }(\cR) -d$. Since the stability condition does
not change under these two kinds of transformations, the
`effective' number of stability parameters of a quiver sheaf
$\cR=(\cE,\phi)$ is $2N(\cR)-2$, where $N(\cR)$ is the (finite) number
of vertices $v\in Q_0$ with $\cE_v\neq 0$. 
From the point of view of the vortex equations \eqref{eq:vortex-equation}, 
the first type of transformations, $\sigma'=c\sigma$, $\tau'=c\tau$,
corresponds to a redefinition of the sections $\phi'=c^{1/2}\phi$
(note that the stability condition is invariant under this
transformation), while the second type corresponds to the constraint
$\deg_{\sigma,\tau }(\cR)=0$ in (i). 
\item[(iii)] 
As usual with stability criteria, in Definition
\ref{def:stab-Q-sheaf}, to check $(\sigma,\tau)$-stability of a
$Q$-sheaf $\cR$, it suffices to check $\mu_{\sigma ,\tau}(\cR')<\mu_{\sigma
,\tau }(\cR)$ for the proper $Q$-subsheaves
$\cR'\subset\cR$ such that $\cE_v'\subset\cE_v$ is saturated, i.e.
such that the quotient $\cE_v/\cE'_v$ is torsion-free, for each $v\in
Q_0$.  
\end{enumerate}}\end{remarks}

\section{Hitchin--Kobayashi correspondence}
\label{sec:HK-Q-bundles}

In this section we will prove a Hitchin--Kobayashi correspondence
between the twisted quiver vortex equations and the stability
condition for holomorphic twisted quiver bundles:

\begin{theorem}\label{thm:HK-Q-bundles}
Let $\sigma$ and $\tau$ be collections of real numbers $\sigma_v$ and
$\tau_v$, respectively, 
with $\sigma_v>0$, for each $v\in Q_0$. Let $\cR=(\cE ,\phi)$ be a
holomorphic $M$-twisted $Q$-bundle such that $\deg_{\sigma, \tau }(\cR)=0$. 
Then $\cR$ is $(\sigma ,\tau)$-polystable if and only if it admits
a hermitian metric $H$ satisfying the quiver
$(\sigma,\tau)$-vortex equations \eqref{eq:vortex-equation}. This
hermitian metric $H$ is unique up to an automorphism of the
$Q$-bundle, i.e. up to a multiplication by a constant
$\lambda_j>0$ for each $(\sigma ,\tau)$-stable summand $\cR_j$ of
$\cR=\cR_1\oplus\cdots\oplus\cR_l$. 
\end{theorem}

\begin{remark}
\label{remark:HK-Q-bundles}{\rm
This theorem generalizes previous theorems, mainly 
Donaldson--Uhlenbeck--Yau theorem \cite{D1,D2,UY}, the 
Hitchin--Kobayashi correspondence for Higgs bundles \cite{H,S},
holomorphic triples and chains \cite{AG1,BG}, twisted holomorphic
triples \cite{BGK2}, etc. 
It should be mentioned that Theorem \ref{thm:HK-Q-bundles} does not
follow from the general theorems proved in \cite{Ba,M} for the
following two reasons. First, the symplectic form $\sum_{v\in
Q_0}\sigma_v\omega_v+\omega_{\RRR}$ on $\AAA\times\Omega^0$
(cf. \S \ref{sub:momentmap}) has been deformed by the parameters
$\sigma$ whenever $\sigma_v\neq\sigma_{v'}$ for some $v,v'\in Q_0$; 
as a matter of fact, the vortex equations \eqref{eq:vortex-equation}
depend on new parameters even for holomorphic triples or 
chains \cite{AG1,BG}, hence generalizing their Hitchin--Kobayashi
correspondences (in the case of a holomorphic pair $(\cE,\phi)$, consisting
of a holomorphic vector bundle $\cE$ and a holomorphic section
$\phi\in H^0(X,\cE)$, as considered in \cite{B}, which can be
understood as a holomorphic triple 
$\phi:\cO_X\to\cE$, the new parameter can actually be absorbed in
$\phi$, so no new parameters are really present). Second, the
twisting bundles $M_a$, for $a\in Q_1$, are not considered in \cite{Ba,M}. 
Our method of proof combines the moment map techniques developed in 
\cite{B,D2,S,UY} for bundles with a proof of a similar
correspondence for quiver modules in \cite[\S 6]{K}.
}\end{remark}

\bfsubsection{Preliminaries and general notation}
\label{sub:notation}

Throughout Section \ref{sec:HK-Q-bundles}, $\cR=(\cE,\phi)$ is a fixed holomorphic
($M$-twisted) $Q$-bundle with $\deg_{\sigma, \tau }(\cR)=0$. 
To prove Theorem \ref{thm:HK-Q-bundles}, 
we can assume that $Q=(Q_0,Q_1)$ is a finite quiver, with $\cE_v\neq
0$, for $v\in Q_0$, and $\phi_a\neq 0$, for $a\in Q_1$ (if this is not
the case, we remove the vertices $v$ with $\cE_v=0$, and the arrows
$a$ with $\phi_a=0$, see Remark \ref{remark:finite-quiver}).
The technical details of the proof largely simplify by introducing
the following notation. 
Unless otherwise stated, $v,v',\ldots$ (resp. $a,a',\ldots$) stand for 
elements of $Q_0$ (resp. $Q_1$), while sums, direct sums and products in
$v,v',\ldots$ (resp. $a,a',\ldots$) are over elements of $Q_0$
(resp. $Q_1$).  
Thus, the condition $\deg_{\sigma, \tau }(\cR)=0$ is equivalent to
$\sum_v\sigma_v\deg(\cE_v)=\sum_v\tau_v\rk(\cE_v)$. 
Let 
\begin{equation}
\label{eq:Eoplus}
\cE=\oplus_v\cE_v; 
\end{equation}
a vector $u$ in the fibre $\cE_x$ over $x\in X$, is a collection
vectors $u_v$ in the fibre $\cE_{v,x}$ over, for each $v\in Q_0$. Let
$\dbar_{\cE_v}:\Omega^0(\cE_v)\to \Omega^{0,1}(\cE_v)$ be the
$\dbar$-operator of the holomorphic vector bundle $\cE_v$, and let 
\begin{equation}
\label{eq:dbarEoplus}
\dbar_{\cE}=\oplus_v\dbar_{\cE_v}
\end{equation}
be the induced $\dbar$-operator on $\cE$. A hermitian metric $H_v$ on
$\cE_v$ defines a 
unique Chern connection $A_{H_v}$ compatible with the holomorphic structure
$\dbar_{\cE_v}$; the corresponding covariant derivative is 
$\du_{H_v}=\partial_{H_v}+\dbar_{\cE_v}$, where $\partial_{H_v}:\Omega^0(\cE_v)
\to\Omega^{1,0}(\cE_v)$ is its $(1,0)$-part. Thus, given
$u\in\Omega^{i,j}(\cE)$, $\dbar_{\cE}(u)\in\Omega^{i,j+1}(\cE)
=\oplus_v\Omega^{i,j+1}(\cE_v)$ is the collection of $\cE_v$-valued
$(i,j+1)$-forms $(\dbar_{\cE}(u))_v =\dbar_{\cE_v}(u_v)$, for each
$v\in Q_0$. 

\bfsubsubsection{Metrics and associated bundles}
\label{subsub:ass-bundles}

Let $Met_v$ be the space of hermitian metrics on $\cE_v$. 
A hermitian metric $(\cdot,\cdot)_{H_v}$ on $\cE_v$ is determined by a
smooth morphism $H_v:\cE_v\to\overline{\cE_v^*}$, by
$(u_v,u'_v)_{H_v}=H_v(u_v)(u'_v)$, with $u_v,u'_v$ in the same fibre of
$\cE_v$. The {\em right} action of the complex gauge group $\GGG^c_v$ on
$Met_v$ is given, by means of this correspondence, by $Met_v\times\GGG^c_v
\to Met_v$, $(H_v,g_v)\mapsto H_v \circ g_v$. 
Let $S_v(H_v)$ be the space of $H_v$-selfadjoint smooth endomorphisms
of $\cE_v$, for each $H_v\in Met_v$. 
We choose a fixed hermitian metric $K_v\in Met$ such that the
hermitian metric $\det(K_v)$ induced by $K_v$ on the determinant
bundle $\det(\cE_v)$ satisfies $\imag\Lambda F_{\det(K_v)}=
\deg(\cE_v)$, for each $v\in Q_0$ (such hermitian metric $K_v$ exists
by Hodge theory). 
Any other metric on $\cE_v$ is given by $H_v=K_v
e^{s_v}$ for some $s_v\in S_v$, or equivalently, by
$(u_v,u'_v)_{H_v}=(e^{s_v}u_v,u'_v)_{K_v}$, where $S_v=S_v(K_v)$.
Let $Met$ be the space of hermitian metrics on $\cE$ such that the
direct sum $\cE=\oplus_v\cE_v$ is orthogonal. A metric $H\in Met$ is given by a
collection of metrics $H_v\in Met_v$, by $(u,u')_H=
\sum_v(u_v,u'_v)_{H_v}$. Let $S(H)=\oplus_v S_v(H_v)$, for
each $H\in Met$, and $S=S(K)=\oplus_v S_v$.  
A vector $s\in S(H)$ is given by a collection of vectors $s_v\in
S_v(H_v)$, for each $v\in Q_0$, while a metric $H\in Met$ is given
by $H=K e^{s}$ for some $s\in S$, i.e. $H_v=K_v e^{s_v}$. 
The (fibrewise) norm on $\cE_v$
(resp. $\cE$) corresponding to $H_v$ (resp.
$H$), is given by $|u_v|_{H_v}=(u_v,u_v)_{H_v}^{1/2}$
(resp. $|u|_{H}=(u,u)_{H}^{1/2}$). The corresponding $L^2$-metric and
$L^2$-norm on the space of sections of $\cE_v$ (resp. $\cE$), is
defined by 
$$
(u_v,u'_v)_{L^2,H_v}=\int_X (u_v,u'_v)_{H_v}, \quad
\| u_v\|_{L^2,H_v}=(u_v,u_v)_{L^2,H_v}^{1/2},
\quad {\rm ~~for~} u_v,u'_v\in \Omega^0(\cE_v),
$$
(resp. $(u,u')_{L^2,H}=\sum_v(u_v,u'_v)_{L^2,H_v}$, 
$\| u\|_{L^2,H}=(u,u)_{L^2,H}^{1/2}$). 
The $L^p$-norm on the space of sections of $\cE$, given by
$$
\| u\|_{L^p,H}=\left(\int_X | u |^p_{H}\right)^{\frac{1}{p}}
\quad {\rm ~~for~} u\in \Omega^0(\cE),
$$
will also be useful.  These metrics
and norms induce canonical metrics on the associated bundles, which
will be denoted with the same symbols. For instance, $H_v\in Met_v$
(resp. $H\in Met$) induces an $L^p$-norm $\|\cdot\|_{L^p,H_v}$ on
$S_v(H_v)$ (resp. $\|\cdot\|_{L^p,H}$ on $S(H)$).
To simplify the notation, we set $(u_v,u'_{v})=(u_v,u'_{v})_{K_v},  
|u_v|=|u_v|_{K_v}, (u,u')=(u,u')_{K}, |u|=|u|_{K}$; 
and $(u_v,u'_{v})_{L^2}=(u_v,u'_{v})_{L^2,K_v}, 
\|u_v\|_{L^2}=\|u_v\|_{L^2,K_v}, (u,u')_{L^2}=(u,u')_{L^2,K}$,
$\|u\|_{L^p}=\|u\|_{L^p,K}$. 

The morphisms $\phi_a:\cE_{ta}\otimes M_a\to\cE_{ha}$ induce a section
$\phi=\oplus_{a} \phi_a$ of the {\em representation bundle}, defined
as the smooth vector bundle over $X$ 
$$
\RRR=\bigoplus_{a} \Hom(\cE_{ta}\otimes M_a,\cE_{ha}). 
$$
A metric $H\in Met$ induces another metric $H_a$ on
each term $\Hom(\cE_{ta}\otimes M_a,\cE_{ha})$ of $\RRR$, by $(\phi_a,
\phi'_a)_{H_a}=\tr(\phi_a\circ\phi_a^{\prime * H_a})$ for $\phi_a,\phi'_a$ 
in the same fibre of $\Hom(\cE_{ta},\cE_{ha})$, where $\phi_a^{\prime
*H_a}:\cE_{ha}\to\cE_{ta}\otimes M_a$ is defined as in \S
\ref{sub:vortex-equation}. Thus, $H$ defines a hermitian metric on
$\RRR$, which we shall also denote $H$, by
$
(\phi,\phi')_{H}=\sum_a (\phi_a,\phi'_a)_{H_a}, 
$
where $\phi,\phi'$ are in a fibre of $\RRR$. The corresponding
fibrewise norm $|\cdot|_H$ is given by $|\phi|_{H}=(\phi,\phi)_{H}^{1/2}$. 
By integrating the hermitian metric over $X$, $(\cdot,\cdot)_{H_a}$
and $(\cdot,\cdot)_{H}$ induce $L^2$-inner products
$(\cdot,\cdot)_{H_a,L^2}$ and $(\cdot,\cdot)_{H,L^2}$ on
$\Omega^0(\cE_{ta}\otimes M_a,\cE_{ha})$ and $\Omega^0=\Omega^0 (\RRR)$
respectively, given by 
$
(\phi_a,\phi_a')_{H_a,L^2}=\int_X (\phi_a,\phi_a)_{H_a}, 
{\rm ~for~} \phi_a,\phi_a'\in\Omega^0 (\cE_{ta}\otimes M_a,\cE_{ha}), 
$ 
and $(\phi,\phi')_{H,L^2}=\sum_a (\phi_a,\phi_a)_{L^2, H_a}, 
{\rm ~for~} \phi,\phi'\in\Omega^0$, 
with associated $L^2$-norms $\|\cdot\|_{H_a,L^2}$, $\|\cdot\|_{H,L^2}$
given by $\|\phi_a\|_{L^2,H} =(\phi_a,\phi_a)_{L^2,H}^{1/2}$ and
$\|\phi\|_{L^2,H} =(\phi,\phi)_{L^2,H}^{1/2}$. 
We set $(\phi,\phi')=(\phi,\phi')_{K}$, 
$|\phi|=|\phi|_{K}$, for each $\phi,\phi'$ in the same fibre
of $\RRR$; and $(\phi,\phi')_{L^2}=(\phi,\phi')_{L^2, K}$, 
$\|\phi\|_{L^2}=\|\phi\|_{L^2,K}$, for each $\phi,\phi'$ smooth
sections of $\RRR$.

\bfsubsubsection{The vortex equations}
\label{subsub:vortex-equations}

Composition of two endomorphisms 
$s, s'\in S$ is defined by $(s\circ s')_v=s_v\circ
s'_v$ for $v\in Q_0$. The identity endomorphism $\id$ of
$\cE$ is given by $\id_v=\id_{\cE_v}$.
Given a vector bundle $F$ on $X$, we define the endomorphisms $\sigma,\tau:
F\otimes\End(\cE)\to F\otimes\End(\cE)$, where $\End(\cE)$ is
the bundle of smooth endomorphisms of $\cE$, by fibrewise multiplication,
i.e. $(\sigma\cdot(f\otimes s))_v= f\otimes\sigma_v s_v$ and
$(\tau\cdot(f\otimes s))_v= f\otimes\tau_v s_v$, for $f\in F$ and 
$s\in\End(\cE)$ in the fibres over the same point $x\in X$. For
instance, if $s\in S$, then $(\sigma\cdot\dbar_{\cE}(s))_v
=\sigma_v\dbar_{\cE_v}(s_v)$.
Given $H\in Met$ and sections $\phi,\phi'$ of $\RRR$, we define the 
endomorphisms $\phi\circ\phi^{\prime * H},\phi^{* H} \circ\phi',
[\phi,\phi^{\prime *H}]\in \Omega^0(\End(\cE))$, using \S \ref{sub:vortex-equation}, 
by 
$$
(\phi\circ\phi^{\prime *H})_v=\sum_{v\in h^{-1}(a)} \phi_a\circ\phi^{\prime *H_a}_a,
\quad 
(\phi^{*H}\circ\phi')_v=\sum_{v\in t^{-1}(a)}\phi^{* H_a}_a\circ\phi'_a,
$$
$$
[\phi,\phi^{\prime *H}]
=\phi\circ\phi^{\prime *H}-\phi^{*H}\circ\phi'.
$$
Note that $[\phi,\phi^{*H}]\in S(H)$. The quiver vortex equations
\eqref{eq:vortex-equation} can now be written in a compact form 
\begin{equation}\label{eq:bvortex-equation}
\sigma\cdot\imag\Lambda F_{H}+[\phi,\phi^{*H}]=\tau\cdot\id , \quad
{\rm for~} H\in Met .
\end{equation}
Given $s\in S$ and $\phi\in\Omega^0=\Omega^0(\RRR)$, $s\circ\phi,\phi\circ
s,[s,\phi],[\phi, s]\in\Omega^0$ are defined by  
$$
(s\circ\phi)_a=s_{ha}\circ\phi_a, ~
(\phi\circ s)_a=\phi_a\circ (s_{ta}\otimes \id_{M_a}), ~
[s,\phi]=s\circ\phi-\phi\circ s, ~
[\phi, s]=\phi\circ s-s\circ\phi.
$$

\bfsubsubsection{The trace and trace free parts of the vortex
equations} 
\label{subsub:trace-part-vortex-eq}

The trace map is defined by $\tr:\End(\cE)\to\CC$,
$s\mapsto\tr(s)=\sum_{v}\tr(s_v)$. Let $S^0(H)$ be the space of
`$\sigma$-trace free' $H$-selfadjoint endomorphisms $s\in S(H)$, i.e.
such that $\tr(\sigma\cdot s)=0$, or more explicitly, 
$\sum_v\sigma_v\tr(s_v)=0$, for each $H\in Met$; let
$S^0=S^0(K)\subset S$.
Let $Met^0$ be the space of metrics $H=K e^s$ with $s\in
S^0$. The metrics $H\in Met^0$ satisfy the trace part of equation
\eqref{eq:bvortex-equation}, i.e. 
\begin{equation}\label{eq:Met-circ}
\tr(\sigma\cdot\imag \Lambda F_{H})=\tr(\tau\cdot \id).
\end{equation}
To prove this, let $H=K e^{s}\in Met$ with $s\in S$. Then 
$\det(H_v)=\det(K_v) e^{\tr s_v}$ so $\tr F_{H_v}=F_{\det(H_v)}=
F_{\det(K_v)}+\dbar\partial\tr s_v=\tr F_{K_v}
+\dbar\partial\tr s_v$ (since the operators induced by
$\dbar_{\det(\cE_v)}$ and $\partial_{\det(K_v)}$ on the trivial bundle
of endomorphisms of $\det(\cE_v)$ are $\dbar$ and
$\partial$, resp.). Adding for all $v$, $ \tr(\sigma\cdot\imag\Lambda
F_{H}) =  \tr(\sigma\cdot\imag\Lambda
F_{K})+\imag\Lambda\dbar\partial\tr(\sigma\cdot s)$, 
where $\tr(\imag\Lambda
F_{K_v})=\deg(\cE_v)$ by construction (cf. \S
\ref{subsub:ass-bundles}), so $\tr(\sigma\cdot\imag\Lambda 
F_{K})=\sum_v\sigma_v\deg(\cE_v)=\sum_v\tau_v\rk(\cE_v)
=\tr(\tau\cdot\id)$. Thus, 
\begin{equation}\label{eq:Met-circ-prime}
\tr(\sigma\cdot\imag\Lambda F_{H}-\tau\cdot\id)
=\imag\Lambda\dbar\partial\tr(\sigma\cdot s), 
\end{equation}
which is zero if $s\in S^0$. This proves \eqref{eq:Met-circ}. 
Therefore, a metric $H=K e^{s}\in Met^0$ satisfies the quiver
$(\sigma,\tau)$-vortex equations \eqref{eq:bvortex-equation} if and
only if it satisfies the trace free part, i.e.  
$$
p^0_H\left(\sigma\cdot\imag\Lambda F_{H}+[\phi,\phi^{*H}]
-\tau\cdot\id\right)=0,
$$
where $p^0_H:S(H)\to S(H)$ is the $H$-orthogonal projection onto
$S^0(H)$. 

\bfsubsubsection{Sobolev spaces} 
\label{subsub:sobolev}

Following \cite{UY,S,B}, given a smooth vector bundle $E$, and any
integers $k,p\geq 0$, $L^p_k\Omega^{i,j}(E)$ is the Sobolev space of
sections of class $L^p_k$, i.e. $E$-valued $(i,j)$-forms whose
derivatives of order $\leq k$ have finite $L^p$-norm. 
Throughout the proof of Theorem \ref{thm:HK-Q-bundles}, we fix an
even integer $p>\dim_\RR(X)=2n$.
Note that there is a compact embedding of 
$L^p_2\Omega^{i,j}(E)$ into the space of continuous $E$-valued
$(i,j)$-forms on $X$, for $p>2n$. This embedding will be
used in \S \ref{subsub:const-involving-met}. 
Particularly important are the collection $L^p_2 S=\oplus_v L^p_2
S_v$ of Sobolev spaces $L^p_2 S_v$ of $K_v$-selfadjoint endomorphisms
of $\cE_v$ of class $L^p_2$; the collection $Met^p_2\cong\prod_v
Met^p_{2,v}$ of Sobolev metrics, with 
$$
Met^p_{2,v}=\{ K_v e^{s_v} |s_v\in L^p_2 S_v\}, \quad
{\rm ~for~each~} v\in Q_0;
$$
the subspace $L^p_2 S^0\subset L^p_2 S$ of sections $s\in
L^p_2 S$ such that $\tr(\sigma\cdot s)=0$ almost everywhere in $X$; and 
$$
Met^{p,0}_2=\{ Ke^{s} |s_v\in L^p_2 S^0\}\subset Met^p_2.
$$
Given $H=Ke^s\in Met^p_2$, with $s\in L^p_2 S$, we define the
$H$-adjoint of $\phi$, generalizing the case where $s_v$ is smooth,
i.e. $\phi^{*H}=e^{-s}\circ \phi^{*K}\circ e^s$.
Similar generalizations apply to the
other constructions in \S\S \ref{subsub:vortex-equations},
\ref{subsub:trace-part-vortex-eq}, to define $L^p_2 S_v(H_v)$ and
$L^p_2 S(H)=\oplus_v L^p_2 S_v(H)$, as well as the subspace 
$L^p_2 S^0(H)\subset L^p_2 S(H)$, for each $H\in Met^p_2$. 
If $H_v=K_v e^{s_v}\in Met^p_{2,v}$ with 
$s_v\in L^p_2 S_v$, we define the connection $A_{H_v}$, with $L^p_1$
coefficients, and its curvature $F_{H_v}\in L^p\Omega^{1,1}(\End(\cE_v))$, 
with $L^p$ coefficients, 
generalizing the case where $s_v$ is smooth: 
\begin{equation}\label{eq:Sobolev-connections}
\du_{H_v}:=\du_{K_v}+e^{-s_v}\partial_{K_v}(e^{s_v}), \quad
F_{H_v}=F_{K_v}+\dbar_{\cE_v}(e^{-s_v}\partial_{K_v}(e^{s_v})), 
\end{equation}
(where $\du_{H_v}$ is the covariant derivative associated to the
connection $A_{H_v}$). 

\bfsubsubsection{The degree of a saturated subsheaf}
\label{subsub:deg-subsheaf}
 
A saturated coherent subsheaf $\cF'$ of a holomorphic vector bundle
$\cF$ on $X$ (i.e., a coherent subsheaf with $\cF/\cF'$ 
torsion-free), is reflexive, hence a vector subbundle
outside of codimension 2. Given a hermitian metric $H$ on $\cF$, the
$H$-orthogonal projection $\pi'$ from $\cF$ onto $\cF'$, defined
outside codimension 2, is an $L^2_1$-section of the bundle
of endomorphisms of $\cF$, so $\beta=\dbar_\cF(\pi')$
is of class $L^2$, where $\dbar_{\cF}$ is the
$\dbar$-operator of $\cF$. The degree of $\cF'$ is 
$$ 
\deg(\cF')=\frac{1}{\Vol(X)}\left(\int_X\tr(\pi'\imag\Lambda F_{H})
-\|\beta\|^2_{L^2,H}\right) ,
$$
(cf. \cite{UY,S,B}).

\bfsubsubsection{Some constructions involving hermitian matrices}
\label{subsub:const-involving-met}

The following definitions slightly generalize \cite[\S 4]{S}. 
Let $\varphi:\RR\to\RR$ and $\varPhi:\RR\times\RR\to
\RR$ be smooth functions. Given $s\in S$, we define $\varphi(s)\in S$
and linear maps $\varPhi(s):S\to S$ and $\varPhi(s):\Omega^0(\RRR)\to\Omega^0(\RRR$) 
(we denote the last two maps with the same symbol since there will not be 
possible confusion between them). Actually, we define maps of fibre
bundles $\varPhi:S\to S(\End\cE)$ and $\varPhi:S\to S(\End\RRR)$, for
certain spaces $S(\End\cE)$ and $S(\End\RRR)$, which we first define. 
Let $S(\End\cE)=\oplus_v S(\End\cE_v)$, where $S(\End\cE_v)$ is
the space of smooth sections of the bundle $\End(\End\cE_v)$ which are
selfadjoint w.r.t. the metric induced by $K_v$. Let $\End\RRR$ be
the endomorphism bundle of the vector bundle $\RRR$; 
$S(\End\RRR)$ is the space of smooth sections of $\End\RRR$ which are
selfadjoint w.r.t. the metric induced by $K_v$ and $q_a$.
We define $\varphi(s_v)\in S_v$ for $s_v\in S_v$ and a linear map
$\varPhi:S_v\to S(\End\cE_v)$ as follows. Let $s_v\in S_v$. If $x\in
X$, let $(u_{v,i})$ be an orthonormal basis of $\cE_{v,x}$ (w.r.t. $K_v$),
with dual basis $(u^{v,i})$, such that $s_v=\sum_{i} \lambda_{v,i} u_{v,i}\otimes
u^{v,i}$. Furthermore, let $(m^{a,k})$ be the dual of an orthonormal
basis of $M_{a,x}$ (w.r.t. $q_a$). The value of $\varphi(s_v)\in S_v$
at the point $x\in X$ is defined as in \cite[\S 4]{S}, by 
\begin{equation}\label{eq:const-involving-met-1}
\varphi(s_v)(x):=\sum_{i} \varphi(\lambda_{v,i}) u_{v,i}\otimes
u^{v,i}.
\end{equation}
We define $\varphi(s)\in S$, for $s\in S$, by $\varphi(s)_v:=\varphi(s_v)$. 
Given $f_v\in S_v$ with $f_{v}(x)=\sum_{i,j} f_{v,ij}
u_{v,i}\otimes u^{v,j}$, the value of $\varPhi(s_v)f_v\in S_v$ at the point
$x\in X$ is 
\begin{equation}\label{eq:const-involving-met-2}
\varPhi(s_v)f_v(x):=\sum_{i,j} \varPhi(\lambda_{v,i},\lambda_{v,j})
f_{v,ij} u_{v,i}\otimes u^{v,j}, 
\end{equation}
and we define $\varPhi:S\to S(\End\cE)$ and $\varPhi:S\to S(\End\RRR)$
as follows. Let $s\in S$. First, if $f\in S$, $(\varPhi(s)f)_v 
:=\varPhi(s_v)f_v$. Second, given a section $\phi$ of $\RRR$ such that
the value of $\phi_a:\cE_{ta}\otimes M_a\to \cE_{ha}$ at $x\in X$ is
$\phi_{a}(x)=\sum_{i,j,k} \phi_{a,ijk}(x) u_{ha,j}\otimes
u^{ta,i}\otimes m^{a,k}$ for each $a\in Q_1$, the value of
$\varPhi(s)\phi\in\Omega^0(\RRR)$ at $x\in X$ is 
\begin{equation}\label{eq:const-involving-met-3}
(\varPhi(s)\phi(x))_a:=\sum_{i,j,k}\varPhi(\lambda_{ha,j},\lambda_{ta,i}) 
\phi_{a,ijk}(x) u_{ha,j}\otimes u^{ta,i}\otimes m^{a,k}, 
\quad {\rm for~each~} a\in Q_1. 
\end{equation}
Note that if $\varPhi$ is given by $\varPhi(x,y)=\varphi_1(x)\varphi_2(y)$
for certain functions $\varphi_1,\varphi_2:\RR\to\RR$, then 
$(\varPhi(s) \phi)_a=\varphi_1(s_{ha})\circ\phi_a\circ(\varphi_2(s_{ha})
\otimes\id_{M_a})$, that is, 
\begin{equation}
  \label{eq:const-involving-met-3-bis}
  \varPhi(s)\phi=\varphi_1(s)\circ\phi\circ\varphi_2(s).
\end{equation} 
Finally, given a smooth function $\varphi:\RR\to\RR$, we define
$\du\varphi:\RR\times\RR\to\RR$ as in \cite[\S 4]{S}:
$$
\du\varphi(x,y)=\frac{\varphi(y)-\varphi(x)}{y-x}, 
{\rm ~if~} x\neq y, {\rm ~and~} \du\varphi(x,y)=\varphi'(x) 
{\rm ~if~} x=y. 
$$
Thus, 
\begin{equation}\label{eq:const-involving-met-4}
\dbar_{\cE}(\varphi(s))=\du\varphi(s)(\dbar_{\cE}(s))
\quad {\rm for~} s\in S. 
\end{equation}
The following lemma will be especially important in the proof of Lemma
\ref{lemma:LimitingEndomorphism-quivers}. Given a number $b$,
$L^2_{k,b}S\subset L^p_k S$ is the closed subset of sections 
$s\in L^2_kS$ such that $|s|\leq b$ a.e. in $X$;
$L^2_{0,b}S(\End\RRR)$ is similarly defined.

\begin{lemma}
\label{lemma:const-involving-met}
\begin{enumerate}
\item[(i)] 
$\varphi:S\to S$ extends to a continous map $\varphi:L^2_{0,b}S\to
L^2_{0,b'}S$ for some $b'$.
\item[(ii)] 
$\varphi:S\to S$ extends to a map $\varphi:L^2_{1,b}S\to
L^q_{1,b'}S$ for some $b'$, for $q\leq 2$, which is continuous for
$q<2$. Formula \eqref{eq:const-involving-met-4} holds in this context. 
\item[(iii)] 
$\varPhi:S\to S(\End\cE)$ extends to a map
$\varPhi:L^2_{0,b}S\to \Hom(L^2\Omega^0(\End\cE), L^q\Omega^0(\End\cE))$ 
for $q\leq 2$, which is continous in the norm operator topology for
$q<2$. 
\item[(iv)] 
$\varPhi:S\to S(\End\RRR)$ extends to a continuous map
$\varphi:L^2_{0,b}S\to L^2_{0,b'}S(\End\RRR)$ for some $b'$.
\item[(v)] 
The previous maps extend to smooth maps $\varphi:L^p_2S\to L^p_2S$, 
$\varPhi:L^p_2 S\to L^p_2 S(\End\cE)$ and $\varPhi:L^p_2 S\to L^p_2
S(\End\RRR)$ between Banach spaces of Sobolev sections. Formulas
\eqref{eq:const-involving-met-1}-\eqref{eq:const-involving-met-4}
hold everywhere in $X$. 
\end{enumerate}
\end{lemma}

\proof 
This follows as in \cite{B,S}. For (v), $p> 2n$, so there is a
compact embedding $L^p_2\subset C^0$. 
\qed

\subsection{Existence of special metric implies polystability}
\label{sub:vortexeq-polystability}

Let $H$ be a hermitian metric on $\cR$ satisfying the
quiver $(\sigma ,\tau)$-vortex equations. To prove that $\cR$ is $(\sigma
,\tau)$-polystable, we can assume that it is indecomposable ---then
we have to prove that it is actually $(\sigma ,\tau)$-stable.
Let $\cR'=(\cE',\phi')\subset\cR$ be a proper $Q$-subsheaf. We can
assume that $\cE'_v\subset\cE_v$ is saturated for each $v\in Q_0$ (cf.
Remark \ref{remark:stability}(iii)).  Let $\pi'_v$ be the
$H_v$-orthogonal projection from $\cE_v$ onto $\cE'_v$, defined
outside codimension 2, $\pi''_v=\id-\pi'_v$, and
$\beta_v=\dbar_\cE(\pi'_v)$. The collections of 
sections $\pi'_v, \pi''_v, \beta_v$ define elements $\pi',\pi''\in
L^2_1\Omega^{0}(\End\cE)$, $\beta\in L^2\Omega^{0,1}(\End\cE)$, respectively. 
Taking the $L^2$-product with $\pi'$ in \eqref{eq:bvortex-equation},
$$
(\sigma\cdot\imag\Lambda F_H,\pi')_{L^2,H}+([\phi,\phi^{*H}],\pi')_{L^2,H}
=(\tau\cdot\id,\pi')_{L^2,H}. 
$$
We now evaluate the three terms of this equation. The first term in
the left hand side is 
$$
(\sigma\cdot\imag\Lambda F_H,\pi')_{L^2,H}
=\sum_v\sigma_v(\imag\Lambda F_{H_v},\pi'_v)_{L^2,H_v}
=\Vol(X)\sum_v\sigma_v\deg(\cE_v)+\sum_v\sigma_v\|\beta_v\|^2_{L^2,H_v}
$$
(cf. \S \ref{subsub:deg-subsheaf}). 
Let $\phi'=\pi'\circ\phi\circ\pi'$, $\phi''=\pi''\circ\phi\circ\pi'$,
$\phi^\perp=\pi'\circ\phi\circ\pi''$. Then 
$\phi=\phi'\circ\pi' +\phi^\perp\circ\pi''+\phi''\circ\pi''$ outside
of codimension 2, for $\cR'\subset\cR$. Thus, 
$[\pi',\phi]=\phi^\perp\circ\pi''$, and the second term is
$$
([\phi,\phi^{*H}],\pi')_{L^2,H}=(\phi,[\pi',\phi])_{L^2,H}
                            =(\phi,\phi^\perp)_{L^2,H}=\|\phi^\perp\|^2_{L^2,H} .
$$
Finally, the right hand side is
$$
(\tau\cdot\id,\pi')_{L^2,H}=\int_X\sum_v\tau_v\tr(\pi'_v)=\Vol(X)\sum_v\tau_v\rk(\cE'_v),
$$ 
(since $\tr(\pi'_v)=\rk(\cE'_v)$ outside of codimension 2). Therefore
$$
\Vol(X)\deg_{\sigma ,\tau }(\cR')
= -\sum_{v\in Q_0}\sigma_v\|\beta_v\|^2_{L^2,H_v}
  -\sum_{a\in Q_1}\|\phi^\perp_a\|^2_{L^2,H_a}.
$$
The indecomposability of $\cR$ implies that either $\beta_v\neq 0$
for some $v\in Q_0$ or $\phi^\perp_a\neq 0$ for some $a\in Q_1$; thus,
$\deg_{\sigma ,\tau }(\cR')<0$, so $\mu_{\sigma ,\tau }(\cR')
<0=\mu_{\sigma ,\tau }(\cR)$, hence $\cR$ is
$(\sigma,\tau)$-stable. 
\qed

\bfsubsection{The modified Donaldson lagrangian}
\label{sub:modified-Donaldson-lagrangian}

To define the modified Donaldson Lagrangian, we first recall the
definition of the Donaldson lagrangian (cf. \cite[\S 5]{S}). Let
$\Psi:\RR\times\RR\to\RR$ be given by  
\begin{equation}\label{eq:Psi-quivers} 
\Psi (x,y)=\frac{e^{y-x}-(y-x)-1}{(y-x)^2}. 
\end{equation}
The {\em Donaldson lagrangian} $M_{D,v}=M_D(K_v,\cdot):
Met^p_{2,v}\to\RR$ is given by 
\begin{align*}
M_{D,v}(H_v) 
= (\imag\Lambda F_{K_v},s_v)_{L^2} 
+ (\Psi(s_v)(\dbar_{\cE_v}s_v),\dbar_{\cE_v}s_v)_{L^2},
{\rm ~~for~} H_v=K_v e^{s_v}\in Met^p_{2,v}, ~ s_v\in L^p_2 S_v.
\end{align*} 
The Donaldson lagrangian $M_{D,v}=M_D(K_v,\cdot)$ is additive in the
sense that 
\begin{equation}
\label{eq:M_D-additive}
M_{D,v}(K_v,H_v)+M_{D,v}(H_v,J_v)=M_{D,v}(K_v,J_v), 
\quad {\rm for~} H_v,J_v\in Met^p_2. 
\end{equation}
Another important property is that the Lie derivative of
$M_{D,v}$ at $H_v\in Met^p_2$, in the direction of $s_v\in L^p_2
S_v(H_v)$, is given by the moment map (cf. \S \ref{sub:momentmap}), i.e. 
\begin{equation}
\label{eq:M_D}
\frac{\du}{\du\varepsilon}M_{D,v}(H_v e^{\varepsilon s_v})\big|_{\varepsilon=0}
=(\imag\Lambda F_{H_v},s_v)_{L^2,H_v}, 
\quad{\rm with ~} H_v\in Met^p_2, s_v\in L^p_2 S_v(H_v).
\end{equation}
Higher order Lie derivatives can be easily evaluated. Thus, from
\eqref{eq:Sobolev-connections}, 
\begin{equation}
\label{eq:M_D-2-prelim}
\frac{\du}{\du\varepsilon} F_{H_v e^{\varepsilon s_v}}=
\dbar_{\cE_v}\partial_{H_v e^{\varepsilon s_v}}s_v, \quad {\rm for~
each~} H_v\in Met^p_2 {\rm ~and~} s_v\in L^p_2 S(H_v)
\end{equation} 
so the second order Lie derivative is 
\begin{equation}
\label{eq:M_D-2}
\frac{\du^2}{\du\varepsilon^2}
M_{D,v}(H_v e^{\varepsilon s_v})\big|_{\varepsilon=0}
=(\imag\Lambda\dbar_{\cE_v}\partial_{H_v}s_v,s_v)_{L^2,H_v}
=\|\dbar_{\cE_v}s_v\|_{L^2,H_v} 
\end{equation}
(the second equality is obtained by integrating $\tr(s_v\imag\Lambda
\dbar_{\cE_v}\partial_{H_v}s_v)=\imag\Lambda\dbar
\tr(s_v\partial_{H_v}s_v)+|\dbar_{\cE_v}s_v|_{H_v}^2$ over $X$,where 
$|\dbar_{\cE_v}s_v|_{H_v}^2= -\imag\Lambda\tr(\dbar_{\cE_v}s_v
\wedge\partial_{H_v}s_v)$ by the \kah\ identities, 
and $\int_X \Lambda\dbar\tr(s_v\partial_{H_v}s_v)= \int_X
\dbar\tr(s_v\partial_{H_v}(s_v))\wedge\omega^{n-1}/(n-1)!=0$ by Stokes
theorem --- cf. e.g. \cite[Lemma 3.1(b) and the proof Proposition
5.1]{S}). 

\begin{definition}\label{def:donaldson-lagrangian}
The {\em modified Donaldson lagrangian} $M_{\sigma,\tau}=M_{\sigma,\tau}(K,\cdot):
Met^p_2\to\RR$ is 
$$
M_{\sigma,\tau}(H)=\sum_{v} \sigma_v M_{D,v}(H_v)
 +\|\phi\|^2_{L^2,H}-\|\phi\|^2_{L^2,K}
-(s,\tau\cdot\id)_{L^2}, \quad {\rm for~} H=Ke^s\in Met^p_2, ~ s\in L^p_2 S.
$$
\end{definition}

Using the constructions of \S \ref{subsub:const-involving-met}, the
modified Donaldson lagrangian can be expressed in terms of the  
functions $\Psi,\psi:\RR\times\RR\to\RR$, with $\Psi$ given by
\eqref{eq:Psi-quivers} and $\psi$ defined by 
\begin{equation}\label{eq:psi-quivers} 
\psi(x,y)=e^{x-y}. 
\end{equation}
In the following, we use the notation
$(\cdot,\cdot)_{L^2}=(\cdot,\cdot)_{L^2,K}$, $\|\cdot\|_{L^2}
=\|\cdot\|_{L^2,K}$, as defined in \S \ref{subsub:ass-bundles}.

\begin{lemma}\label{lemma:bpsis}
If $H=K e^{s}\in Met^p_2$, with $s\in L^p_2 S$, then 
$$
M_{\sigma,\tau}(H) = (\sigma\cdot\imag\Lambda F_{K},s)_{L^2}
+(\sigma\cdot\Psi(s)(\dbar_{\cE}s),\dbar_{\cE}s)_{L^2} 
+ (\psi(s)\phi,\phi)_{L^2}-\|\phi\|^2_{L^2}
-(\tau\cdot\id,s)_{L^2}.
$$
\end{lemma}
 
\proof 
The first two terms follow from the definitions of $M_{D,v}$ and
$M_{\sigma,\tau}$. To obtain the third term, we note that 
$\phi_a^{*H_{a}}=(e^{-s_{ta}}\otimes\id_{M_a})\circ \phi_a^{*K_a}\circ
e^{s_{ha}}$ and $(\psi(s)\phi)_a=e^{s_{ha}} \circ\phi_a\circ
(e^{-s_{ta}}\otimes\id_{M_a})$ (cf. \eqref{eq:const-involving-met-3-bis}), 
so $|\phi_a |^2_{H_a}=\tr(\phi_a\circ\phi_a^{*H_a})
=\tr(e^{s_{ha}}\circ\phi_a\circ(e^{-s_{ta}}\otimes\id_{M_a})\circ\phi_a^{*K_a})
=\tr((\psi(s)\phi)_a\circ\phi_a^{*K_a})=((\psi(s)\phi)_a,\phi_a)_{K_a}$.
The last two terms follow directly from the definition of $M_{\sigma,\tau}$. 
\qed

\subsection{Minima of $\bmM_{\bsigma,\btau}$, the main estimate, and
  the vortex equations} 
\label{sub:minima}

Let $m_{\sigma,\tau}:Met^p_2\to L^p\Omega^0(\End\cE)$ be defined by 
\begin{equation}\label{eq:def-m}
m_{\sigma,\tau}(H)=\sigma\cdot\imag\Lambda F_{H}+[\phi,\phi^{*H}]
-\tau\cdot\id , \quad {\rm for~} H=Ke^s\in Met^p_2, ~ s\in L^p_2 S.
\end{equation}
Thus, $m_{\sigma,\tau}(H)\in L^pS(H)$ for each $H\in Met^p_2$, and
actually $m_{\sigma,\tau}(H)\in L^pS^0(H)$ if $H\in Met^{p,0}_2$, by
\eqref{eq:Met-circ}. 
Let $B>\|m_{\sigma,\tau} (K)\|^p_{L^p}$ be a positive real number. 
We are interested in the minima of $M_{\sigma,\tau}$ in the closed subset of
$Met^{p,0}_2$ defined by 
$$
\Met^{p,0}_{2,B}:=\{H\in\Met^{p,0}_2|\,\|m_{\sigma,\tau} (H)\|^p_{L^p,H}
                       \leq B\} 
$$
(the restriction to this subset will be necessary to apply Lemma
\ref{boundedfunctions} below). 

\begin{proposition}
\label{prop:minima}
If $\cR$ is simple, i.e. its only endomorphisms are multiples of the
identity, and $H\in Met^{p,0}_{2,B}$ minimises $M_{\sigma,\tau}$ on
$\Met^{p,0}_{2,B}$, then $m_{\sigma,\tau}(H)=0$. 
\end{proposition}

The minima are thus the solutions of the vortex equations.
To prove this, we need a lemma about the first and second order Lie 
derivaties of $M_{\sigma,\tau}$. Given $H\in Met^p_2$, $L_H:L^p_2
S(H)\to L^p S(H)$ is defined by 
\begin{equation}\label{eq:L-1}
L_{H}(s)=\frac{\du}{\du\varepsilon}
m_{\sigma,\tau} (He^{\varepsilon s})\big|_{\varepsilon=0}, 
\quad {\rm for~each~} s\in L^p_2 S(H). 
\end{equation}
Since $\phi^{*H_\varepsilon}= e^{-\varepsilon s}\phi^{*H} e^{\varepsilon
s}$, with $H_\varepsilon=H e^{\varepsilon s}$, we have 
\begin{equation}
\label{eq:dphi*}
\frac{\du}{\du\varepsilon}\phi^{*H_\varepsilon}\big|_{\varepsilon=0}
=[s,\phi]^{*H}, 
\end{equation}
so $\frac{\du}{\du\varepsilon}[\phi,\phi^{*H_\varepsilon}]\big|_{\varepsilon=0}
=[\phi,[s,\phi]^{*H}]$.  
Together with \eqref{eq:M_D-2-prelim}, this implies that 
\begin{equation}\label{eq:L-2}
L_H(s)=\sigma\cdot\imag\dbar_{\cE}\partial_{H}s + [\phi,[s,\phi]^{*H}].
\end{equation}

\begin{lemma}
\label{lemma:derivatives-M}
\begin{enumerate}
\item[(i)] 
${\displaystyle 
M_{\sigma,\tau}(K,H)+M_{\sigma,\tau}(H,J)=M_{\sigma,\tau}(K,J)}$,
for $H,J\in Met^p_2$; 
\item[(ii)]
${\displaystyle 
\frac{\du}{\du\varepsilon} M_{\sigma,\tau} (He^{\varepsilon
s})\big|_{\varepsilon=0} =(m_{\sigma,\tau}(H),s)_{L^2,H}}$, for each
$H\in Met^p_2$ and $s\in L^p_2 S(H)$; 
\item[(iii)]
${\displaystyle \frac{\du^2}{\du\varepsilon^2}
M_{\sigma,\tau} (He^{\varepsilon s})\big|_{\varepsilon=0}
=(L_H(s),s)_{L^2,H}
=\sum_v\sigma_v\|\dbar_{\cE_v}s_v\|^2_{L^2,H_v}+\|[s,\phi]\|^2_{L^2,H}}$,
for each $H\in Met^p_2$ and $s\in L^p_2 S(H)$.
\end{enumerate}
\end{lemma}

\proof
Part (i) follows immediately from \eqref{eq:M_D-additive} and
$(Ke^s)e^{s'}=Ke^{s+s'}$. 
To prove (ii) and (iii), let $H_{\varepsilon}=He^{\varepsilon s}$, for
$\varepsilon\in\RR$. From \eqref{eq:dphi*} we get
$\frac{\du}{\du\varepsilon}|\phi|^2_{H_\varepsilon}\big|_{\varepsilon=0}
= \tr\left(\phi \frac{\du}{\du\varepsilon}
\phi^{*H_\varepsilon}\big|_{\varepsilon=0}\right) =\tr(\phi
[s,\phi]^{*H})=([\phi,\phi^{*H}],s)_H$, which together with
\eqref{eq:M_D}, proves (ii) (the last term in \eqref{eq:def-m} is
trivially obtained). The first equality in (iii) follows from (ii),
the $H_\varepsilon$-selfadjointness of $s$ (since 
$s^{*H_\varepsilon}=e^{-\varepsilon s} s^{* H}e^{\varepsilon
s}=e^{-\varepsilon s} s e^{\varepsilon s}=s$), and \eqref{eq:L-1}: 
$$
\frac{\du^2}{\du\varepsilon^2} M_{\sigma,\tau} (H_\varepsilon)\big|_{\varepsilon=0}
=\frac{\du}{\du\varepsilon} (m_{\sigma,\tau}(H_\varepsilon),s)_{L^2,H_\varepsilon}\big|_{\varepsilon=0}
=\!\int_X\tr\left(\frac{\du}{\du\varepsilon}m_{\sigma,\tau}(H_\varepsilon)\big|_{\varepsilon=0} s\right)
=\!\int_X \tr(L_H(s)s), 
$$
which equals $(L_H(s), s)_{L^2,H}$. To prove the second equality in
(iii), we first notice that if $\phi'$ is a smooth section of $\RRR$,
then $(s,\phi'\circ\phi^{*H})_H=(s\circ\phi,\phi')_H$ and
$(s,\phi^{*H}\circ\phi')_H=(\phi\circ s, \phi')_H$, so
$(s,[\phi',\phi^{*H}])_H=([s,\phi],\phi')_H$. The second equality in
(iii) is now obtained using \eqref{eq:L-2}, \eqref{eq:M_D-2-prelim} and
taking $\phi'=[s,\phi]$ in the previous formula. 
\qed

\begin{proof}[Proof of Proposition \ref{prop:minima}.]
We start proving that if $\cR$ is simple and $H\in Met^{p,0}_2$, then
the restriction of $L_H$ to $L^p_2 S^0(H)$, which we also denote by $L_H:
L^p_2 S^0(H)\to L^p S^0(H)$, is surjective. 
To do this, we only have to show that $L_H$ is a Fredholm operator of
index zero and that it has no kernel. First, for each vertex $v$,
$k_v:L^p_2 S_v(H_v)\to L^p S_v(H_v)$, defined by $k_v=
\imag\Lambda\dbar_{\cE_v}\partial_{H_v}-\imag\Lambda\dbar_{\cE_v}\partial_{K_v}$, 
is obviously a compact operator (cf. \S \ref{subsub:sobolev}), and by
the \kah\ identities, $\imag\Lambda\dbar_{\cE_v}\partial_{K_v}$ acting
on $L^p_2 S$ is the $(1,0)$-laplacian $\Delta'_{K_v}=\partial_{K_v}\partial_{K_v}^*
+\partial_{K_v}^*\partial_{K_v}$, which is elliptic and selfadjoint,
hence Fredholm, and has index zero. Now, $L_H$ equals 
$\sum_v \sigma_v\imag\Lambda\dbar_{\cE}\partial_{H_v}$, up to a
compact operator, so it is also a Fredholm operator of index zero. To
prove that it has no kernel, we notice that if $s\in L^p_2 S^0(H)$
satisfies $L_H(s)=0$, then $(s,L_H(s))_{L^2,H}=0$, so Lemma
\ref{lemma:derivatives-M}(iii) implies $\dbar_{\cE_v}s_v=0$ and
$[s,\phi]=0$; 
i.e. $s$ is actually an endomorphism of $\cR$, so
$s_v=c\id_{\cE_v}$, for certain constant $c$. Since
$\tr(\sigma\cdot s)=0$, the constant is $c=0$, so $s_v=0$. 

Let $H$ minimise $M_{\sigma,\tau}$ in $Met^{p,0}_{2,B}$. 
To prove that $m_{\sigma,\tau}(H)=0$, we assume the contrary. 
Since $L_H:L^p_2 S^0(H)\to L^p S^0(H)$ is surjective, and
$m_{\sigma,\tau}(H)\in 
S^0(H)$ is not zero, there exists a non-zero 
$s\in L^p_2 S^0(H)$ with 
$L_H(s)= -m_{\sigma,\tau}(H)$. We shall consider the values of
$M_{\sigma,\tau}$ along the path $H_{\varepsilon}=H e^{\varepsilon s}
\in Met^{p,0}_2$ for small $|\varepsilon|$. First, 
$$
\frac{\du}{\du\varepsilon} |m_{\sigma,\tau}(H_{\varepsilon})|^2_{H_{\varepsilon}}\big|_{\varepsilon=0}
= \frac{\du}{\du\varepsilon} \tr(m_{\sigma,\tau}(H_{\varepsilon})^2)\big|_{\varepsilon=0}
=2(m_{\sigma,\tau}(H),L_H(s))_H=-2|m_{\sigma,\tau}(H)|^2_H, 
$$
(cf. \eqref{eq:L-1}), and since $p$ is even, 
$$
\frac{\du}{\du\varepsilon}\| m_{\sigma,\tau}(H_{\varepsilon})\|^p_{L^p,H_{\varepsilon}}\big|_{\varepsilon=0} 
=\frac{p}{2} \int_X |m_{\sigma,\tau}(H)|^{p-2}_{H}
\frac{\du}{\du\varepsilon} |m_{\sigma,\tau}(H_{\varepsilon})|^2_{H_{\varepsilon}}\big|_{\varepsilon=0}
=-p\| m_{\sigma,\tau}(H)\|^p_{L^p,H}<0, 
$$
so the path $H_\varepsilon$ is in $Met^{p,0}_{2,B}$ for
small $|\varepsilon|$. Thus, $\frac{\du}{\du\varepsilon}M_{\sigma,\tau}
(H_\varepsilon)\big|_{\varepsilon=0} =0$, as $H$ minimises $M_{\sigma,\tau}$ in
$Met^{p,0}_{2,B}$. Now, Lemma \ref{lemma:derivatives-M}(ii) applied to
$s\in L^p_2 S(H)$ gives 
$$
\frac{\du}{\du\varepsilon} M_{\sigma,\tau}(H_\varepsilon)\big|_{\varepsilon=0}
=(m_{\sigma,\tau}(H),s)_{L^2,H}
=-(L_{H}(s),s)_{L^2,H}.
$$
As in the first paragraph of this proof, if $\cR$ is simple and $s\in
L^p_2 S^0(H)$ satisfies $(s,L_{H}(s))_{L^2,H}=0$, then Lemma
\ref{lemma:derivatives-M}(iii) implies that $s$ is
zero. This contradicts the assumption $m_{\sigma,\tau}(H)\neq 0$. 
\end{proof}

\begin{definition}\label{def:main-estimate}
We say that $M_{\sigma,\tau}$ satisfies the main estimate in $\Met^{p,0}_{2,B}$ if 
there are constants $C_1,C_2>0$, which only depend on $B$, such that  
$\sup |s|\leq C_1 M_{\sigma,\tau}(H)+C_2$, for all $H=K e^{s}
\in\Met^{p,0}_{2,B}$, $s\in L^p_2 S$.
\end{definition}

\begin{proposition}\label{prop:estimate-vortex}
If $\cR$ is simple and $M_{\sigma,\tau}$ satisfies the main estimate in
$\Met^{p,0}_{2,B}$, then there is a hermitian metric on $\cR$
satisfying the $(\sigma,\tau)$-vortex equations. This hermitian metric 
is unique up to multiplication by a positive constant.
\end{proposition}

\proof
This result is proved in exactly the same way as in \cite[\S 3.14]{B},
so here we only sketch the proof.
One first shows that if $M_{\sigma,\tau}(Ke^s)$ is bounded above, then the
Sobolev norms $\| s\|_{L^p_2}$ are bounded. One then takes a minimising
sequence $\{ Ke^{s_j}\}$ for $M_{\sigma,\tau}$, with $s_j\in L^p_2
S^0$; then $\| s_j\|_{L^p_2}$ are uniformly bounded, so after passing
to a subsequence, $\{ s_j\}$ converges weakly in $L^p_2$ to some $s$.
One then sees that $M_{\sigma,\tau}$ is continuous in the weak topology on
$Met^{p,0}_{2,B}$, so $M_{\sigma,\tau}(Ke^{s_j})$ converges to
$M_{\sigma,\tau}(Ke^{s})$. Thus, $H=Ke^s$ minimises $M_{\sigma,\tau}$.
By Proposition \ref{prop:minima}, $m_{\sigma,\tau}(H)=0$, i.e. $H$
satisfies the vortex equations. By elliptic regularity, $H$ is smooth.
The uniqueness of the solution $H$ follows from the convexity of
$M_{\sigma,\tau}$ (cf. Lemma \ref{lemma:derivatives-M}(iii)) and the
simplicity of $\cR$. 
\qed

The proof of Theorem \ref{thm:HK-Q-bundles} is therefore reduced to
show that if $\cR$ is $(\sigma,\tau)$-stable, then $M_{\sigma,\tau}$
satisfies the main estimate in $\Met^{p,0}_{2,B}$ (this is the content of
\S \ref{sub:main-estimate}).

\subsection{Equivalence of $\bmC^{\bf 0}$ and $\bmL^{\bf 1}$ estimates} 
\label{sub:equiv-C0-L1}

The following proposition will be used in \S \ref{sub:main-estimate}. 

\begin{proposition} \label{subprop:equiv-sup-L1} 
There are two constants $C_1,C_2>0$, depending on $B$ and $\sigma$, 
such that for all $H=K e^{s} \in\Met^{p,0}_{2,B}$, $s\in L^p_2 S^0$, 
$\sup |s| \leq C_1 \|s\|_{L^1}+C_2$.
\end{proposition}

\begin{corollary} \label{corollary:main-estimate-L1}
$M_{\sigma,\tau}$ satisfies the main estimate in $\Met^{p,0}_{2,B}$ if and
only if there are constants $C_1,C_2>0$, which only depend on $B$, such that  
$\|s\|_{L^1}\leq C_1 M_{\sigma,\tau}(H)+C_2$, for all $H=K e^{s}
\in\Met^{p,0}_{2,B}$, $s\in L^p_2 S^0$.
\qed
\end{corollary}

Corollary \ref{corollary:main-estimate-L1} is immediate from Proposition 
\ref{subprop:equiv-sup-L1}. To prove Proposition \ref{subprop:equiv-sup-L1}, 
we need three lemmas. The first one is due to Donaldson \cite{D3} (see
also the proof of \cite[Proposition 2.1]{S}). 

\begin{lemma}\label{boundedfunctions}
There exists a smooth function $a:\, [0,\infty)\to [0,\infty )$, with
$a(0)=0$ and $a(x)=x$ for $x>1$, such that the following is true:
For any $\wt{B}\in\RR$, there is a constant $C(\wt{B})$ such that if 
$f$ is a positive bounded function on $X$ and $\Delta f\leq b$, 
where $b$ is a function in $L^p(X)~(p>n)$ with $\|b\|_{L^p}\leq \wt{B}$,
then $\sup |f|\leq C(\wt{B})a(\| f\|_{L^1})$. Furthermore, if $\Delta f\leq
0$, then $\Delta f=0$. 
\qed
\end{lemma}

\begin{lemma}
\label{lemma:k-leq-h}
If $s\in L^p_2 S$ and $H=K e^{s}\in\Met^p_2$, then   
$([\phi,\phi^{*H}],s)\geq ([\phi,\phi^{*K}],s)$.
\end{lemma}

\proof  
The function $f(\varepsilon)=([\phi,\phi^{*H_\varepsilon}],s)$ for
$\varepsilon\in\RR$, where $H_\varepsilon=Ke^{\varepsilon s}$, is
increasing, as ${\du }f(\varepsilon)/{\du\varepsilon}
=|[s,\phi]|^2_{H_{\varepsilon}} \geq 0$ (cf. \eqref{eq:dphi*}).
Now,$f(0)=([\phi,\phi^{*K}],s)$, $f(1)=([\phi,\phi^{*H}],s)$, so we are done.
\qed 

\begin{lemma}\label{lemma:delta-m} 
If $H=K e^{s}\in\Met^p_{2}$, with $s\in L^p_2 S$, then
$$
(m_{\sigma,\tau}(H)-m_{\sigma,\tau}(K),s)\geq \frac{1}{2} ~
|\sigma^{1/2}\cdot s| ~ \Delta |\sigma^{1/2}\cdot s|,
$$
where $\sigma^{1/2}\cdot s\in L^p_2 S$ is of course defined by
$(\sigma^{1/2}\cdot s)_v=\sigma^{1/2}_v s_v$, for $v\in Q_0$. 
\end{lemma}

\proof
This lemma, and its proof, are similar to (but not completely
immediate from) \cite[Proposition 3.7.1]{B}. First, Lemma
\ref{lemma:k-leq-h} and \eqref{eq:Sobolev-connections} imply
\begin{equation}
\label{eq:delta-m-1}
(m_{\sigma,\tau}(H)-m_{\sigma,\tau}(K),s)\geq 
\imag \Lambda (\sigma\cdot F_H - \sigma\cdot F_K,s)
=\imag\Lambda (\sigma\cdot \dbar_\cE(e^{-s}\partial_K e^s),s),
\end{equation}
where 
\begin{equation}
\label{eq:delta-m-2}
(\sigma\cdot \dbar_\cE(e^{-s}\partial_K e^s),s)
=\dbar(\sigma\cdot e^{-s}\partial_K e^s,s)
        +(\sigma\cdot e^{-s}\partial_K e^s,\partial_K s)
\end{equation}
(for $A_K$ is the Chern connection corresponding to the metric $K$).
To make some local calculations, we choose a local $K_v$-orthogonal 
basis $\{ u_{v,i} \}$ of eigenvectors of $s_v$, for each vertex $v$,
with corresponding eigenvalues $\{\lambda_{v,i}\}$, and let $\{
u^{v,i} \}$ be the corresponding dual basis; thus,
$$
\label{eq:delta-m-3}
s_v=\sum_i \lambda_{v,i} u_{v,i} \otimes u^{v,i}. 
$$
As in \cite[(3.36)]{B}, a local calculation gives
$(e^{-s_v}\partial_{K_v} e^{s_v},s_v)=\frac{1}{2}\partial |s_v|^2$; 
multiplying by $\sigma_v$ and adding for $v\in Q_0$, we get 
$(\sigma\cdot e^{-s}\partial_{K} e^{s},s)=\frac{1}{2}\partial |s'|^2$, 
where $s'=\sigma^{1/2}\cdot s$. 
Thus, 
\begin{equation}
\label{eq:delta-m-6}
\dbar(\sigma\cdot  e^{-s}\partial_K e^s,s)=\frac{1}{2}\dbar\partial |s'|^2
=|s'|\dbar\partial|s'|+\dbar|s'|\wedge\partial|s'|
\end{equation}
From \eqref{eq:delta-m-1}, \eqref{eq:delta-m-2}, \eqref{eq:delta-m-6} and the
equality $\Delta=2\imag\Lambda\dbar\partial$ for the action of the laplacian
on $0$-forms in a \kah\ manifold, we get
$$
(m_{\sigma,\tau}(H)-m_{\sigma,\tau}(K),s)\geq \frac{1}{2} |s'| \Delta |s'|
        + \imag \Lambda (\dbar |s| \wedge \partial |s'|)
        + \imag \Lambda (\sigma\cdot e^{-s}\partial_{K} e^{s},\partial_K s).
$$
In the proof of \cite[Proposition 3.7.1]{B}, there are several local
calculations which, although there they are only used for the section
$s\in L^p_2 S$ defining the metric $H=K e^s$, are actually valid for
any $K$-selfadjoint section, in particular for $s'\in L^p_2 S$. Thus,
\cite[(3.42)]{B} applied to $s_v$ is 
$$
\imag\Lambda(e^{-s_v}\partial_{K_v} e^{s_v},\sigma\cdot\partial_{K_v} s_v)
\geq \sum_i \imag\Lambda
(\partial\lambda_{v,i}\wedge\dbar\lambda_{v,i}), 
$$
and multiplying by $\sigma_v$ and adding for $v\in Q_0$, we get
\begin{equation}
\label{eq:delta-m-7}
\imag\Lambda(\sigma\cdot e^{-s}\partial_{K} e^{s},\sigma\cdot\partial_{K} s)
\geq \sum_{v,i} \imag\Lambda (\partial\lambda'_{v,i}\wedge\dbar\lambda'_{v,i}), 
\end{equation}
where $\lambda'_{v,i}:=\sigma_v^{1/2} \lambda_{v,i}$ are the eigenvalues of
$s'_v=\sigma_v^{1/2} s_v$; similarly, \cite[(3.43)]{B} applied to $s'$ is
\begin{equation}
\label{eq:delta-m-8}
\sum_{v,i} \imag\Lambda (\partial\lambda'_{v,i}\wedge\dbar\lambda'_{v,i})
\geq \imag\Lambda (\partial|s'|\wedge\dbar|s'|)
=-\imag\Lambda (\dbar|s'|\wedge\partial|s'|).
\end{equation}
From \eqref{eq:delta-m-6}, \eqref{eq:delta-m-7}, \eqref{eq:delta-m-8},
we obtain 
$(m_{\sigma,\tau}(H)-m_{\sigma,\tau}(K),s)\geq \frac{1}{2}
|s'| \Delta |s'|$. 
\qed

\begin{proof}[Proof of Proposition \ref{subprop:equiv-sup-L1}.]
Let $\sigma_{\min}=\min\{\sigma_v |v\in Q_0\}$, 
$\sigma_{\max}=\max\{\sigma_v |v\in Q_0\}$. Given $H=K e^s\in Met^{p,0}_{2,B}$, 
with $s\in L^p_2 S^0$, let $f=|\sigma^{1/2}\cdot s|$ and $b=
2\sigma^{-1/2}_{\min} (|m_{\sigma,\tau}(H)|+|m_{\sigma,\tau}(K)|)$.
We now verify that $f$ and $b$ satisfy the hypotheses of Lemma
\ref{boundedfunctions}, for a certain $\wt{B}$ which only depends
on $B$. First, $\| b\|_{L^p}\leq 2\sigma^{-1/2}_{\min} 
(\|m_{\sigma,\tau}(H)\|_{L^p} +\|m_{\sigma,\tau}(K)\|_{L^p})
\leq \wt{B}:=2\sigma^{-1/2}_{\min} 2 B^{1/p}$. Second, we prove that 
\begin{equation}\label{eq:equiv-sup-L1}
\Delta f \leq b.
\end{equation}
At the points where $f$ does not vanish, $f^{-1}\leq
\sigma^{-1/2}_{\min} |s|^{-1}$, so Lemma \ref{lemma:delta-m} gives 
$$
\Delta f \leq 2\sigma^{-1/2}_{\min} |s|^{-1} (m_{\sigma,\tau}(H)-m_{\sigma,\tau}(K),s)
\leq 2\sigma^{-1/2}_{\min} |m_{\sigma,\tau}(H)-m_{\sigma,\tau}(K)|\leq b, 
$$
while to consider the points where $f$ vanishes, we just take into
account that $\Delta f=0$ almost everywhere (a.e.) in
$f^{-1}(0)\subset X$, and that $b\geq
0$ by its definition, so \eqref{eq:equiv-sup-L1} actually holds a.e.
in $X$. The hypotheses of Lemma \ref{boundedfunctions} are thus
satisfied, so there exists a constant $C(B)>0$ such that $\sup f\leq
C(B) a(\| f\|_{L^1})$, with with $a:[0,\infty)\ra [0,\infty)$ as in Lemma
\ref{boundedfunctions}. This estimate can also be written as
$\sup f \leq C_1 \| f \|_{L^1} + C_2$, where $C_1, C_2>0$ only
depend on $B$. Now, $|s|\leq \sigma^{-1/2}_{\min} f$ and $f\leq
\sigma^{1/2}_{\max}|s|$, so 
$$ 
\sup|s|  \leq \sigma^{-1/2}_{\min} (C_1 \| f \|_{L^1} + C_2)
  \leq \sigma^{-1/2}_{\min} (C_1 \sigma^{1/2}_{\max} \| s \|_{L^1} + C_2)
$$
The estimate is obtained by redefining the constants $C_1,C_2$.
\end{proof}

\bfsubsection{Stability implies the main estimate}
\label{sub:main-estimate}

The following proposition, together with Proposition 
\ref{prop:estimate-vortex}, are the key ingredients to complete the
proof of Theorem \ref{thm:HK-Q-bundles} (cf. Definition
\ref{def:main-estimate} for the main estimate).

\begin{proposition}\label{prop:main-estimate-quivers} 
If $\cR$ is $(\sigma,\tau)$-stable, then $M_{\sigma,\tau}$ satisfies
the main estimate in $\Met^{p,0}_{2,B}$. 
\end{proposition}

To prove this, we need some preliminaries (Lemmas
\ref{sublemma:sequence-quivers}-\ref{prop:weak-subsheaves}). 
Let $\{ C_j\}_{j=1}^\infty$ be a sequence of constants with
${\displaystyle \lim_{j\to\infty}C_j=\infty}$. 

\begin{lemma}\label{sublemma:sequence-quivers} 
If $M_{\sigma,\tau}$ does not satisfy the main estimate in
$\Met^{p,0}_{2,B}$, then there is a sequence $\{s_j\}_{j=1}^{\infty}$ in
$L^p_2S^0$ with $Ke^{s_j}\in Met^{p,0}_{2,B}$ (which we can assume to
be smooth), such that 
\begin{enumerate}
\item[(i)] ${\displaystyle \lim_{j\to\infty} \| s_j\|_{L^1}=\infty}$, 
\item[(ii)] ${\displaystyle \| s_j\|_{L^1}\geq C_j M(K e^{s_j})}$.
\end{enumerate}
\end{lemma}

\proof
Let ${b}> \|m_{\sigma,\tau}(K)\|^p_{L^p}$ with ${b}<B$, so
$\Met^{p,0}_{2,{b}}\subset\Met^{p,0}_{2,B}$. Thus, if 
$M_{\sigma,\tau}$ does not satisfy the main estimate in
$\Met^{p,0}_{2,B}$, then it does not satisfy the main estimate in
$\Met^{p,0}_{2,{b}}$ either. 
We shall prove that for any positive constant $C'$, if there are
positive constants $C''$ and $N$ such that $\|s\|_{L^1}\leq C'
M_{\sigma,\tau}(K e^{s})+ C''$ whenever $s\in L^p_2 S^0$ with
$Ke^s\in\Met^{p,0}_{2,{b}}$ and $\| s\|_{L^1}\geq N$, then
$M_{\sigma,\tau}$ satisfies the main estimate in $\Met^{p,0}_{2,{b}}$.
The lemma follows from this claim by choosing a
sequence of constants $\{ N_j\}_{j=1}^\infty$ with $N_j\to\infty$, and
taking $C''_j$ and $s_j\in L^p_2 S^0$ with $Ke^{s_j}\in 
Met^{p,0}_{2,{b}}\subset\Met^{p,0}_{2,B}$, $\| s_j\|_{L^1}\geq N_j$,
and $\|s\|_{L^1}> C_j M_{\sigma,\tau}(K e^{s_j})+ C''_j$. 
Let $C',C'',N$ be such that 
$$
\|s\|_{L^1}\leq C' M_{\sigma,\tau}(K
e^{s})+ C'' \text{~for~} \| s\|_{L^1}\geq N.
$$ 
Let $S_{N}=\{ s\in L^p_2 S^0|Ke^s\in\Met^{p,0}_{2,{b}} \textrm{~and~}
\|s\|_{L^1}\leq N\}$.   
By Proposition \ref{subprop:equiv-sup-L1}, if $s\in S_{N}$, then 
$\sup|s_v|\leq\sup|s|\leq C_1 \|s\|_{L^1}+C_2\leq C_1 N+C_2$ (here
$C_1$ and $C_2$ are not the first elements of the sequence $\{
C_j\}_{j=1}^\infty$ but constants as in Proposition
\ref{subprop:equiv-sup-L1}), so by Lemma \ref{lemma:bpsis},
$M_{\sigma,\tau}$ is bounded below on $S_{N}$, i.e.
$M_{\sigma,\tau}(Ke^s)\geq -\lambda$ for each $s\in 
S_{N}$, for some constant $\lambda>0$. Thus, $\| s\|_{L^1}\leq
C'(M_{\sigma,\tau}(K e^{s})+\lambda)+N$ for each $s\in S_{N}$. 
Replacing $C''$ by $\max\{ C'',C'\lambda+N\}$, we see that 
$\|s\|_{L^1}\leq C' M_{\sigma,\tau}(Ke^s)+C''$, for each $s\in L^p_2
S^0$ with $Ke^s\in\Met^{p,0}_{2,{b}}$. By Corollary
\ref{corollary:main-estimate-L1}, $M_{\sigma,\tau}$ satisfies the main
estimate in $\Met^{p,0}_{2,{b}}$. Finally, since the set of smooth sections is
dense in $L^p_2 S^0$, we can always assume that $s_j$ is smooth (we
made the choice ${b}<B$ so that if $Ke^{s_j}$ is in the boundary
$\|m_{\sigma,\tau} (H)\|^p_{L^p,H}={b}$ of $Met^{p,0}_{2,{b}}$, we can
still replace $s_j$ by a smooth $s'_j$ with $Ke^{s_j'}\in Met^{p,0}_{2,B}$). 
\qed

\begin{lemma}
\label{lemma:LimitingEndomorphism-quivers}
Assume that $M_{\sigma,\tau}$ does not satisfy the main estimate in
$\Met^{p,0}_{2,B}$. Let $\{s_j\}_{j=1}^{\infty}$ be a 
sequence as in Lemma \ref{sublemma:sequence-quivers}, $l_j=\|
s_j\|_{L^1}$, $C(B)=C_1+C_2$, where $C_1, C_2$ are as
in Proposition \ref{subprop:equiv-sup-L1}, and $u_j=s_j/l_j$. Thus, 
$\|u_j\|_{L^1}=1$ and $\sup |u_j|\leq C(B)$.
After going to a subsequence, $u_j\to u_{\infty}$ weakly in
$L^2_1 S^0$, for some nontrivial $u_{\infty}\in L^p_2S^0$ such that
if $\FFF:\RR\times\RR \ra\RR$ is a smooth non-negative function such that
$\FFF(x,y)\leq 1/(x-y)$ whenever $x>y$, and $\FFF_\varepsilon:\RR\times 
\RR\ra\RR$ is a smooth non-negative function with $\FFF_\varepsilon(x,y)=0$
whenever $x-y\leq \varepsilon$, for some fixed $\varepsilon>0$, then 
$$ 
(\sigma\cdot\imag\Lambda F_{K},u_{\infty})_{L^2} 
+ (\sigma\cdot\FFF(u_{\infty})\dbar_{\cE}u_{\infty},
\dbar_{\cE}u_{\infty})_{L^2} \\
+(\FFF_\varepsilon(s)\phi,\phi)_{L^2}
-(\tau\cdot\id,u_{\infty})_{L^2}\leq 0.
$$
\end{lemma}

\proof 
To prove this inequality, we can assume that $\FFF$ and
$\FFF_\varepsilon$ have compact support (for $\sup |u_j|$ are bounded,
by Lemma \ref{subprop:equiv-sup-L1}, and the definitions of
$\FFF(s)\dbar_\cE u_\infty$ and $\FFF_\varepsilon(s)\phi$
only depend on the values of $\FFF$ and $\FFF_\varepsilon$ at the
pairs $(\lambda_i,\lambda_j)$ of eigenvalues, as seen in \S
\ref{subsub:const-involving-met}). Now, if $\FFF$ and
$\FFF_\varepsilon$ have compact support then, for large enough $l$, 
$$
\FFF(x,y)\leq l\Psi(lx,ly), \quad \FFF_\varepsilon(x,y)\leq l^{-1} \psi(lx,ly),
$$
where $\Psi$ and $\psi$ are defined as in \eqref{eq:Psi-quivers} and
\eqref{eq:psi-quivers} (cf. the proof of \cite[Proposition 3.9.1]{B}).
Since $l_j\to\infty$, from these inequalities we obtain that for large
enough $j$, 
$$
(\FFF(u_{j,v})\dbar_{\cE}u_{j,v},\dbar_{\cE}u_{j,v})_{L^2}\leq l
(\Psi(l_{j,v}u_{j,v})\dbar_{\cE}u_{j,v}, \dbar_{\cE}u_{j,v})_{L^2}, 
$$
$$ 
(\FFF_\varepsilon(u_j)\phi,\phi)_{L^2}\leq l^{-1} (\psi(l_ju_j)\phi,\phi)_{L^2},
$$
so Lemma \ref{sublemma:sequence-quivers}(iii) applied to 
$s_i=l_ju_j$, together with lemma \ref{lemma:bpsis}, give an upper
bound 
\begin{align*}
    \frac{1}{C_j}+\frac{\|\phi\|^2_{L^2}}{l_j} &\geq
    l_j^{-1} M_{\sigma,\tau}(Ke^{l_ju_j}) + l_j^{-1} \|\phi\|^2_{L^2}
    \geq (\sigma\cdot\imag\Lambda F_K,u_j)_{L^2} \\
     & +(\sigma\cdot\FFF(u_{j})\dbar_{\cE}u_{j},\dbar_{\cE}u_{j})_{L^2}
    +(\FFF_{\varepsilon}(u_j)\phi,\phi)_{L^2} - (\tau\cdot\id,u_j)_{L^2}.
\end{align*}
As in the proof of \cite[Proposition 3.9.1]{B}, one can use this upper
bound to show that the sequence $\{u_j\}_{j=1}^\infty$
is bounded in $L^2_1$. Thus, after going to a subsequence, $u_j\to
u_\infty$ in $L^2_1$, for some $u_\infty\in L^2_1 S$ with
$\|u_\infty\|_{L^1}=1$, so $u_\infty$ is non-trivial. 

We now prove the estimate for $u_\infty$. First, since $\sup|u_j|\leq
b:=C(B)$, $u_j\to u_\infty$ in $L^2_{0,b}$; applying Lemma
\ref{lemma:const-involving-met}(iii)
, one can show (as in the
proof of \cite[Lemma 5.4]{S}) 
that $(\sigma\cdot\imag\Lambda F_K,u_j)_{L^2} + 
(\sigma\cdot\FFF(u_{j})\dbar_{\cE}u_{j},\dbar_{\cE}u_{j})_{L^2}$
approaches $(\sigma\cdot\imag\Lambda F_K,u_\infty)_{L^2} + (\sigma\cdot\FFF
(u_{\infty})\dbar_{\cE}u_{\infty},\dbar_{\cE}u_{\infty})_{L^2}$ as
$j\to\infty$. Second, since $L^2_1\subset L^2$ is a compact embedding
and actuallly $u_j\in L^2_{1,b}S\subset L^2_{0,b}S$, applying Lemma
\ref{lemma:const-involving-met}(iv) 
(as in the proof of \cite[Proposition 3.9.1]{B}), $\FFF_{\varepsilon}:L^2_{0,b}S \to
L^2_{0,b'}S(\End\RRR)$, $u\mapsto \FFF_{\varepsilon}(u)$, is
continuous on $L^2_{0,b}S$, so $\lim_{j\to \infty} \FFF_{\varepsilon}(u_j) =
\FFF_{\varepsilon}(u_\infty)$. Since $\sup|u_j|$ are bounded, this implies that
$(\FFF_{\varepsilon}(u_j)\phi,\phi)_{L^2}$ converges to
$(\FFF_{\varepsilon}(u_\infty)\phi,\phi)_{L^2}$ as $j\to\infty$.
Finally, it is clear that $(\tau\cdot\id,u_j)_{L^2}\to(\tau\cdot\id,
u_\infty)_{L^2}$ as $j\to\infty$. This completes the proof.
\qed

\begin{lemma}\label{sublemma:eigenvalues-quiver} 
If $M_{\sigma,\tau}$ does not satisfy the main estimate in
$\Met^{p,0}_{2,B}$, and $u_{\infty}\in L^p_2S^0$ is an in Lemma
\ref{lemma:LimitingEndomorphism-quivers}, then the following happens: 
\begin{enumerate}\item[(i)] 
The eigenvalues of $u_\infty$ are constant almost everywhere. 
\item[(ii)] 
Let the eigenvalues of $u_\infty$ be $\lambda_1,\ldots ,\lambda_r$. If
$\FFF\, :\RR\times\RR\lra\RR$ satisfies $\FFF (\lambda_i
,\lambda_j)=0$ whenever $\lambda_i>\lambda_j$, $1\leq i,j\leq r$, then
$\FFF (u_{\infty} )(\dbar_{\cE} u_\infty )=0$.
\item[(iii)] 
If $\FFF_\varepsilon$ is an in Proposition
\ref{lemma:LimitingEndomorphism-quivers}, then
$\FFF_\varepsilon(u_\infty)\phi=0$. 
\end{enumerate}
\end{lemma}

\proof
Parts (i) and (ii) of are proved as in \cite[appendix]{UY}, \cite[\S\S
6.3.4 and 6.3.5]{S}, or \cite[\S\S 3.9.2 and 3.9.3]{B}, using
Lemma\ref{lemma:const-involving-met}(ii) for part (i) and the
estimate in Lemma \ref{lemma:LimitingEndomorphism-quivers} for part
(ii). Part (iii) is similar to \cite[Lemma 3.9.4]{B}, and again uses
the estimate in Lemma \ref{lemma:LimitingEndomorphism-quivers}.
\qed

We now construct a filtration of quiver subsheaves of $\cR$ using
$L^p_2$-subsystems, as in \cite[\S 3.10]{B}.

\begin{lemma}
\label{prop:weak-subsheaves}
Assume that $M_{\sigma,\tau}$ does not satisfy the main estimate in
$\Met^{p,0}_{2,B}$. Let $u_{\infty}\in L^p_2S^0$ be as in Lemma
\ref{lemma:LimitingEndomorphism-quivers}. 
Let the eigenvalues of $u_\infty$, listed in ascending order, be 
$\lambda_0<\lambda_1<\cdots<\lambda_r$. Since $u_\infty$ is `$\sigma$-trace
free' (cf. \S \ref{subsub:trace-part-vortex-eq}), there are at least
two different eigenvalues, i.e. $r\geq 1$. Let $p_0,\ldots,p_r:\RR\ra\RR$ be
smooth functions such that, for $j<r$, $p_j(x)=1$ if  
$x\leq\lambda_j$, $p_j(x)=0$ if $x\geq\lambda_{j+1}$, and 
$p_r(x)=1$ if $x\leq\lambda_r$. Let $\pi_v:\cE\to\cE_v$ be the
canonical projections (cf. \eqref{eq:Eoplus}) and $\dbar_{\cE}$ be as
in \eqref{eq:dbarEoplus}. The operators $\pi'_r=p_j(u_\infty)$ and
$\pi'_{j,v}=\pi'_j\circ\pi_v$, for $0\leq j\leq r$, satisfy: 
\begin{enumerate}
\item[(i)] $\pi'_j\in L^2_1 S$, $\pi^{\prime 2}_j=\pi'_j=\pi^{\prime
*K}_j$ and  $(1-\pi'_j)\dbar_{\cE}\pi'_j=0$.
\item[(ii)] $(\id-\pi'_{j,ha})\circ\phi_a\circ(\pi'_{j,ta}\otimes\id_{M_a})=0$
for each $v\in Q_0$. 
\item[(iii)] Not all the eigenvalues of $u_{\infty}$ are positive. 
\end{enumerate}
\end{lemma}

\proof
The proof of (i) is as in \cite{S} (right below Lemma 5.6; 
see also \cite[Proposition 3.10.2(i)-(iii)]{B}). 
Part (ii) is similar to, but more involved than, \cite[Proposition
3.10.2(iv)]{B}, so we now give a detailed proof of this part. For each
$j$, let $\varepsilon>0$ be such that $\varepsilon\leq(\lambda_{j+1}-
\lambda_j)/2$, and $\varphi_1,\varphi_2:\RR\to\RR$ be smooth
non-negative functions such that $\varphi_1(x)=0$ if $x\leq
\lambda_{j+1}-\varepsilon/2$ and $\varphi_1(x)=1$ if $x\geq
\lambda_{j+1}$, in the case of $\varphi_1$; 
and $\varphi_2(y)=1$ if $y\leq \lambda_{j}$ and $\varphi_2(y)=0$ if 
$y\geq \lambda_{j}+\varepsilon/2$, in the case of $\varphi_2$. Let
$\FFF_{\varepsilon}:\RR\times\RR\to\RR$ be given by 
$$
\FFF_{\varepsilon}(x,y)=\varphi_1(x)\varphi_2(y). 
$$
If $\FFF_{\varepsilon}(x,y)\neq 0$, then
$x>\lambda_{j+1}-\varepsilon/2$ and $y<\lambda_j+\varepsilon/2$, so
$x-y>\lambda_{j+1}-\lambda_j-\varepsilon\geq\varepsilon$; thus,
$\FFF_{\varepsilon}$ satisfies the hypothesis of Lemma 
\ref{sublemma:eigenvalues-quiver}(iii), so $\FFF_\varepsilon(u_\infty)\phi=0$.
But $\FFF_\varepsilon(u_\infty)\phi=\varphi_1(u_\infty)\circ\phi\circ
\varphi_2(u_\infty)$ (cf. \eqref{eq:const-involving-met-3-bis}), where 
$\varphi_1(u_\infty)=\id-\pi'_j$ and $\varphi_2(u_\infty)=\pi'_j$,
which completes the proof of part (ii). Finally, part (iii) follows
from $\tr(\sigma\cdot u_\infty)=0$ and the non-triviality of $u_\infty$. 
\qed 

\begin{proof}[Proof of Proposition \ref{prop:main-estimate-quivers}.]
Assume that $M_{\sigma,\tau}$ does not satisfy the main estimate in
$Met^{p,0}_{2,B}$. We have to prove that $\cR$ is not $(\sigma,\tau)$-stable. 
By Lemma \ref{prop:weak-subsheaves}(i), the operators
$\pi'_{j,v}$ are weak holomorphic vector subbundles of $\cE_v$, for
$v\in Q_0$ \cite[\S 4]{UY}. 
Applying Uhlenbeck--Yau regularity theorem \cite[\S 7]{UY}, 
they represent reflexive subsheaves $\cE'_{j,v}\subset\cE_v$, and by
Lemma \ref{prop:weak-subsheaves}(ii), the inclusions
$\cE'_{j,v}\subset\cE_v$ are compatible with the morphisms
$\phi_a$, hence define $Q$-subsheaves $\cR'_j=(\cE'_j,\phi'_j)$
of $\cR=(\cE,\phi)$. We thus get a filtration of $Q$-subsheaves 
$$
0\hra\cR'_0\hra\cR'_1\hra\cdots\hra\cR'_r=\cR.
$$
As in \cite[(3.7.2)]{B},
$$
u_{\infty}=\lambda_0\pi'_{0}+\sum_{j=1}^{r}\lambda_j(\pi'_{j}-\pi'_{j-1})
          =\lambda_r\id_{\cE}-\sum_{j=0}^{r-1}(\lambda_{j+1}-\lambda_{j})
          \pi'_{j}, 
$$
so the $v$-component $u_{\infty,v}=u_{\infty}\circ\pi_v$ of $u_\infty$ is
\begin{equation} \label{eq-noestimate1-quivers}
u_{\infty,v}=\lambda_r\id_{\cE_v}-\sum_{j=0}^{r-1}(\lambda_{j+1}-\lambda_{j})
\pi'_{j,v},
\end{equation}
(note that it may happen that $\pi'_{j,v}=\pi'_{j+1,v}$ for some $v$
and $j$). From \eqref{eq:const-involving-met-4} and
$\pi'_{j,v}=p_j(u_{\infty,v})$, $\dbar_{\cE_v}\pi'_{j,v}=\du
p_j(u_{\infty,v}) (\dbar_{\cE_v}u_{\infty,v})$, so 
\begin{equation} \label{eq-noestimate4-quivers}\begin{split}
\sum_{j=0}^{r-1}(\lambda_{j+1}-\lambda_j)|\dbar_{\cE_v} \pi'_{j,v}|^2
        &=\sum_{j=0}^{r-1}(\lambda_{j+1}-\lambda_j)
        ((\du p_j)^2(u_{\infty,v})\dbar_{\cE_v}(u_{\infty,v}),
        \dbar_{\cE_v}(u_{\infty,v})) \\
        &=(\FFF
        (u_{\infty,v})(\dbar_{\cE_v}u_{\infty,v}),\dbar_{\cE_v}u_{\infty,v}), 
\end{split}\end{equation} 
where $\FFF :\,\RR\times\RR\lra\RR$, defined by $\FFF 
=\sum_{j=0}^{l-1}(\lambda_{j+1}-\lambda_j)(\du p_j)^2$, satisfies
the conditions of Lemma \ref{lemma:LimitingEndomorphism-quivers} (cf.
e.g. 
the proof of \cite[Lemma 5.7]{S}). 
We make use of the previous calculations to estimate the number
$$
\chi=\Vol(X)
\left(\lambda_r\deg_{\sigma,\tau}(\cR) -\sum_{j=0}^{r-1}(\lambda_{j+1}-\lambda_j)\deg_{\sigma,\tau}(\cR'_j)\right).
$$
On the one hand, the degree of the subsheaf $\cE'_{j,v}\subset\cE_v$
is given by \eqref{subsub:deg-subsheaf}, 
$$
\Vol(X)\deg(\cE'_{j,v})=(\imag\Lambda F_{K_v},\pi'_{j,v})_{L^2}
-\|\dbar_{\cE_v}\pi'_{j,v}\|^2_{L^2}, 
$$ 
and this formula, together with equations \eqref{eq-noestimate1-quivers} 
and \eqref{eq-noestimate4-quivers}, imply 
\begin{equation*}
\begin{split}
\chi & =\!\!\sum_{v\in Q_0} \sigma_v
\left(\imag\Lambda  F_{K_v},\lambda_r\id_{\cE_v}-\sum_{j=0}^{r-1}(\lambda_{j+1}-\lambda_j)\pi'_{j,v}\right)_{L^2} 
 \!\!\!\! +\sum_{v\in Q_0} \sigma_v\sum_{j=0}^{r-1}(\lambda_{j+1}
  -\lambda_j) \|\dbar_{\cE_v}\pi'_{j,v}\|^2_{L^2} \\ & \quad\quad 
-\sum_{v\in Q_0}\tau_v\Vol(X)\left(\lambda_r\rk(\cE_v)-\sum_{j=0}^{r-1}(\lambda_{j+1}-\lambda_j)\rk(\cE'_{j,v})\right)\\
& 
=(\sigma\cdot \imag\Lambda F_K,u_\infty)_{L^2} +
(\sigma\cdot\FFF(u_{\infty})(\dbar_\cE u_\infty),\dbar_\cE u_\infty)_{L^2}
-(\tau\cdot\id,u_\infty)_{L^2}. 
\end{split}
\end{equation*}
It follows from Lemma \ref{lemma:LimitingEndomorphism-quivers}
(with $\FFF_{\varepsilon}=0$, cf. Lemma \ref{sublemma:eigenvalues-quiver}(iii)), 
that $\chi\leq 0$.
On the other hand, if $\cR$ is $(\sigma,\tau)$-stable, then 
$\mu_{\sigma,\tau}(\cR)>\mu_{\sigma,\tau}(\cR'_j)$, 
for $0\leq j< r$, and since $\sigma\cdot u_\infty\in L^p_2 S^0$ is
trace free, 
$$
\tr(\sigma\cdot u_\infty)=\sum_v \sigma_v \tr(u_\infty\circ\pi_v)
=\lambda_r\sum_{v\in Q_0}\sigma_v\rk(\cE_{v})
-\sum_{j=0}^{r-1}(\lambda_{j+1}-\lambda_j)\sum_{v\in
  Q_0}\sigma_v\rk(\cE'_{j,v})=0, 
$$
so we get 
\begin{equation*}
  \begin{split}
\chi &=\frac{\Vol(X)}{\sum_{v\in Q_0}\sigma_v\rk(\cE_v)}
 \sum_{j=0}^{r-1}(\lambda_{j+1}-\lambda_j)\!\left(\sum_{v\in Q_0}
 \sigma_v\rk(\cE'_{j,v})~\deg_{\sigma,\tau}(\cR) -\!\! 
\sum_{v\in Q_0} \sigma_v\rk(\cE_v)\deg_{\sigma,\tau}(\cR'_j)\right) \\ &
=\Vol(X)\sum_{j=0}^{r-1}(\lambda_{j+1}-\lambda_j)\sum_{v\in Q_0} \sigma_v  \rk(\cE'_{j,v})
(\mu_{\sigma,\tau}(\cR)-\mu_{\sigma,\tau}(\cR'_j)) >0. 
\end{split}
\end{equation*}
Therefore, if $M_{\sigma,\tau}$ does not satisfy the main estimate in
$\Met^{p,0}_{2,B}$, then $\cR$ cannot be $(\sigma,\tau)$-stable. 
\end{proof}

\subsection{Stability implies existence and uniquenes of special
  metric}
\label{sub:stability-vortexeq}

Let $\cR=(\cE,\phi)$ be a $(\sigma,\tau)$-polystable holomorphic
$Q$-bundle on $X$. To prove that it admits a hermitian metric
satisfying the quiver $(\sigma,\tau)$-vortex equations, we can assume that 
$\cR$ is $(\sigma,\tau)$-stable, which in particular implies that
it is simple. The existence and uniqueness of a hermitian metric
satisfying the quiver $(\sigma,\tau)$-vortex equations is now
immediate from Propositions \ref{prop:estimate-vortex} and
\ref{prop:main-estimate-quivers}. 
\qed

Sections \ref{sub:vortexeq-polystability} and
\ref{sub:stability-vortexeq} prove Theorem \ref{thm:HK-Q-bundles}. 

\section{Yang--Mills--Higgs functional and Bogomolov inequality} 

Let $\sigma ,\tau$ be collections of real numbers $\sigma_v,\tau_v$,
with $\sigma_v>0$, for $v\in Q_0$. Given a smooth complex vector
bundle $E$, let $c_1(E)$ and $ch_2(E)$ be its first Chern class and
second Chern character, respectively. By Chern--Weil 
theory, if $A$ is a connection on $E$ then $c_1(E)$ (resp. $ch_2(E)$) is
represented by the closed form $\frac{\imag}{2\pi}\tr(F_A)$
(resp. $-\frac{1}{8\pi^2} \tr(F_A^{~2})$). Define the topologial
invariants of $E$ 
\begin{equation}
\label{eq:C1(E)}
C_1(E)=\int_X c_1(E)\wedge\frac{\omega^{n-1}}{(n-1)!}
      =\frac{1}{2\pi}\int_X \tr(\imag\Lambda F_A)\frac{\omega^n}{n!}
\end{equation}
and
\begin{equation}
\label{eq:Ch2(E)}
Ch_2(E)=\int_X ch_2(E)\wedge\frac{\omega^{n-2}}{(n-2)!}
      =-\frac{1}{8\pi^2}\int_X\tr(F_A^{~2})\wedge\frac{\omega^{n-2}}{(n-2)!}
\end{equation}
(thus, $C_1(E)$ is the degree of $E$, up to a normalisation factor). 
Given a holomorphic vector bundle $\cE$ on $X$, we denote by
$C_1(\cE)$ and $Ch_2(\cE)$ the corresponding topological invariants of
its underlying smooth vector bundle.

\begin{theorem}\label{thm:Bogomolov}
If $\cR=(\cE ,\phi)$ is a $(\sigma,\tau)$-stable holomorphic
$Q$-bundle on $X$, and the $q_q$-selfadjoint endomorphism $\imag
\Lambda F_{q_a}$ of $M_a$ is positive semidefinite, for each $a\in
Q_0$, then
\begin{equation}
  \label{eq:Bogomolov}
  \sum_{v\in v}\tau_v C_1(\cE_v)\geq 2\pi\sum_{v\in Q_0}\sigma_v Ch_2(\cE_v).
\end{equation}
If $C_1(\cE_v)=0$, $Ch_2(\cE_v)=0$ for all $v\in
Q_0$, then the connections $A_{H_v}$ are flat for each $v\in Q_0$, and 
\begin{equation}
  \label{eq:stable-quiver-modules}
    \sum_{a\in h^{-1}(v)}\phi_a\circ\phi^{*H}_a -\sum_{a\in 
      t^{-1}(v)}\phi^{*H}_a\circ\phi_a =\tau_v\id_{E_v}
\end{equation}
for each $v\in Q_0$, where $H$ is a solution of the $M$-twisted quiver
$(\sigma,\tau)$-vortex equations on $\cR$.
\end{theorem}

Thus, quiver bundles can be useful to construct flat connections. 
Note that when $X$ is an algebraic variety,
\eqref{eq:stable-quiver-modules} means that $\cR$ is a family of
$\tau$-stable $Q$-modules parametrized by $X$ (cf. \cite[\S\S 5, 6]{K}). 

This theorem is an immediate consequence of the Hitchin--Kobayashi
correspondence for holomorphic $Q$-bundles and Proposition
\ref{prop:YMH-breaking} below. We shall use the notation introduced in
\S \ref{sub:momentmap}. 

\begin{definition}
\label{def:YMH}
The Yang--Mills--Higgs functional $YMH_{\sigma,\tau}:\AAA\times\Omega^0
\ra\RR$ is defined by
\begin{multline*}
YMH_{\sigma ,\tau}(A,\phi)=\sum_{v\in Q_0} \sigma_v \| F_{A_v}\|^2_{L^2} +
\sum_{a\in Q_1} \| \du_{A_a}\phi_a\|^2_{L^2}\\
+ 2 \sum_{v\in Q_0}\sigma_v^{-1}
\left\|\sum_{a\in h^{-1}(v)}\phi_a\circ\phi^{*H}_a-\sum_{a\in t^{-1}(v)}\phi^{*H}_a\circ\phi_a -\tau_v\id_{E_v}\right\|^2_{L^2}, 
\end{multline*}
where $A_a$ is the connection induced by $A_{ta}$, $A_{q_a}$ and
$A_{ha}$ on the vector bundle $\Hom(E_{ta}\otimes M_a ,E_{ha})$.
\end{definition}

In the following, $\|\cdot\|$ will mean the $L^2$-norm in the
appropiate space of sections. Note that in Theorem \ref{thm:Bogomolov}
it is assumed that $\imag \Lambda F_{q_a}$ is semidefinite positive
for each $a\in Q_0$, so it defines a semidefinite positive sesquilinear
form on $\Omega^0(\Hom(E_{ta} \otimes M_a, E_{ha}))$ by 
$$
(\phi_a,\phi'_a)_{q_a}=\!\!\!\int_X \!\!\!
\tr\left(\phi_a\circ (\id_{E_{ta}}\otimes \imag \Lambda F_{q_a}) \circ  \phi_a^{*H_a}\right), 
\textnormal{ for each } 
\phi_a,\phi'_a\in \Omega^0(\Hom(E_{ta} \otimes M_a, E_{ha})). 
$$
Adding together, we thus get a semidefinite positive sesquilinear form
on $\Omega^0$, defined by
$$
(\phi,\phi')_{\RRR, M}=\sum_{a\in Q_1} (\phi_a,\phi'_a)_{L^2,q_a}, 
\textnormal{  for each  } \phi,\phi'\in \Omega^0.
$$
Thus, $\|\phi\|^2_{\RRR, M}:=(\phi,\phi)_{\RRR, M}\geq 0$ for each
$\phi\in \Omega^0$. 

\begin{proposition}\label{prop:YMH-breaking}
If $(A,\phi)\in\cA'\times\Omega^0$, with $A_v\in\cA^{1,1}_v$
for all $v\in Q_0$, then
\begin{multline*}
YMH_{\sigma,\tau}(A,\phi)
=4\sum_{a\in Q_1} \|\dbar_{A_a}\phi_a\|^2
+4\pi \sum_{v\in Q_0} \tau_v C_1(E_v) 
- 8\pi^2 \sum_{v\in Q_0} \sigma_v Ch_2(E_v) 
- \|\phi\|^2_{\RRR,M}\\
+\sum_{v\in Q_0} \sigma^{-1}_v\left\|\sigma_v\imag\Lambda F_{A_v}
+\sum_{a\in h^{-1}(v)}\phi_a\circ\phi^{*H}_a
-\sum_{a\in t^{-1}(v)}\phi^{*H}_a\circ\phi_a-\tau_v\id_{E_v}\right\|^2.
\end{multline*}
\end{proposition}

\proof
Before giving the proof, we need several preliminaries. First, note
that for any $A_v\in\cA^{1,1}_v$, 
\begin{equation}\label{eq:YMH-1}
\| F_{A_v} \|^2=\|\Lambda F_{A_v}\|^2 -8\pi^2 Ch_2(E_v)
\end{equation}
(cf. e.g. \cite[Theorem 4.2]{B}). Secondly, we notice that the
curvature of $A_a$, for $A\in Q_1$, is given 
by
\begin{equation}\label{eq:YMH-2}
F_{A_a}(\phi_a)=F_{A_{ha}}\circ\phi_a-\phi_a\circ
  (F_{A_{ta}}\otimes\id_{M_a} + \id_{E_{ta}}\otimes F_{q_a})
\end{equation}
where $\phi_a$ is a section of $\Hom(E_{ta},E_{ha})$. Finally, since
the $(0,1)$-parts of the unitary connections $A_{ta}, A_{ha}$ define
holomorphic structures, $A_a$ also defines a holomorphic structure on
the smooth vector bundle $\Hom(E_{ta},E_{ha})$, so it satisfies the
\kah\ identities 
$$
\imag [\Lambda,\partial_{A_a}]=-\dbar^*_{A_a},\quad
\imag [\Lambda,\dbar_{A_a}]=\partial^*_{A_a}.
$$
In particular, the commutator of $\imag\Lambda$ with the curvature $F_{A_a}=
\partial_{A_a}\dbar_{A_a}+\dbar_{A_a}\partial_{A_a}$ is
$\imag[\Lambda,F_{A_a}]=\Delta'_{A_a}-\Delta''_{A_a}$, where
$\Delta'_{A}=\partial_A^*\partial_A +\partial_A\partial_A^*$ and
$\Delta''_{A}=\dbar_A^*\dbar_A +\dbar_A\dbar_A^*$. When acting on
sections $\phi_a$ of $\Hom(E_{ta},E_{ha})$, this simplifies to
$$
\imag\Lambda F_{A_a}\phi_a=\Delta'_{A_a}\phi_a-\Delta''_{A_a}\phi_a.
$$
so that
\begin{equation}\label{eq:YMH-3}
(\imag\Lambda
F_{A_a}\phi_a,\phi_a)_{L^2}=\|\partial_{A_a}\phi_a\|^2-\|\dbar_{A_a}\phi_a\|^2.
\end{equation}
To prove the proposition, we define
$$
U_v(\phi)=\sum_{a\in h^{-1}(v)}\phi_a\circ\phi^{*H}_a -
\sum_{a\in t^{-1}(v)}\phi^{*H}_a\circ\phi_a
$$
for $\phi\in\Omega^0$ and $v\in Q_0$. Then
\begin{multline*}
\sum_{v\in Q_0} \sigma^{-1}_v\|\sigma_v\imag\Lambda F_{A_v}+
U_v(\phi)-\tau_v\id_{E_v}\|^2
= \sum_{v\in Q_0}  \sigma_v \|\Lambda F_{A_v}\|^2 \\
+ \sum_{v\in Q_0}  \sigma^{-1}_v \| U_v(\phi)-\tau_v\id_{E_v}\|^2
+ 2 \sum_{v\in Q_0}  (\imag\Lambda F_{A_v},U_v(\phi))_{L^2}
- 2 \sum_{v\in Q_0} \sigma_v^{-1} (\imag\Lambda F_{A_v},\tau_v\id_{E_v})_{L^2},
\end{multline*}
where \eqref{eq:YMH-2}, \eqref{eq:YMH-3} give
\begin{align*}
\sum_{v\in Q_0}  & (\imag\Lambda F_{A_v},U_v(\phi))_{L^2}
=\sum_{a\in Q_1} (\imag\Lambda F_{A_{ha}}\circ\phi_a 
-\phi_a\circ(\imag\Lambda F_{A_{ta}}\otimes\id_{M_a}),\phi_a)_{L^2}  \\
&=\sum_{a\in Q_1} (\imag\Lambda F_{A_a}\phi_a,\phi_a)_{L^2} -\|\phi\|_{\RRR,M}
=\sum_{a\in Q_1}  \|\partial_{A_a}\phi_a\|^2-\sum_{a\in Q_1} \|\dbar_{A_a}\phi_a\|^2-\|\phi\|_{\RRR,M} .
\end{align*}
The proposition now follows from the previous equation, 
\eqref{eq:YMH-1}, and the definition of $C_1(E_v)$.
\qed

\begin{proof}[Proof of Theorem \ref{thm:Bogomolov}.]
Let $\cR=(\cE,\phi)$ be $(\sigma ,\tau)$-stable, $H$ the
hermitian metric on $\cR$ satisfying the $(\sigma,\tau)$-vortex
equations (cf. Theorem \ref{thm:HK-Q-bundles}), and 
$A\in\AAA$ the corresponding Chern connection. By
Definition \ref{def:YMH}, $\YMH_{\sigma,\tau}(A,\phi)\geq 0$, 
while from Proposition \ref{prop:YMH-breaking}, this is $2\pi
\sum_{v\in Q_0} \tau_v C_1(E_v) - 8\pi^2 \sum_{v\in Q_0} \sigma_v
Ch_2(E_v) -\|\phi\|^2_{\RRR,M}$, as $\dbar_{A_a}\phi_a=0$ for each
$a\in Q_1$. Since we are assuming $\|\phi\|^2_{\RRR,M}\geq 0$, we
obtain \eqref{eq:Bogomolov}. Furthermore, if
$C_1(\cE_v)=Ch_2(\cE_v)=0$ for each $v\in Q_0$, then
$\YMH_{\sigma,\tau}(A,\phi)=-\|\phi\|^2_{\RRR,M}\leq 0$, but this
functional is non-negative by Definition \ref{def:YMH}, so
$\YMH_{\sigma,\tau}(A,\phi)=0$. Thus, $F_{A_v}=0$ and 
we also obtain \eqref{eq:stable-quiver-modules} for each $v\in Q_0$,
again by Definition \ref{def:YMH}. 
\end{proof}

\section{Twisted quiver sheaves and path algebras} \label{sec:pathalgebra}

The category of $M$-twisted $Q$-sheaves is equivalent to the
category of coherent sheaves of right $\euscA$-modules, where $\euscA$
is certain locally free 
$\cO_X$-sheaf associated to $Q$ and $M$ ---the so-called {\em
$M$-twisted path algebra} of $Q$. 
This provides an alternative point of view of twisted quiver sheaves which, in
certain cases, gives a more algebraic understanding of
certain properties of $Q$-sheaves (cf. e.g. \S \ref{sub:tensorprod}
below). In particular, it may be a better point of view to study 
the moduli space problem, which we will not address in this paper.
To fix terminology, a locally free (resp. free, coherent)
$\cO_X$-algebra is a sheaf $\euscS$ of rings which at the same time is a
locally free (resp. free, coherent) $\cO_X$-module. Given such an
$\cO_X$-algebra $\euscS$, a locally free  (resp. free, coherent)
$\euscS$-algebra is a sheaf $\euscA$ of (not necessarily commutative)
rings over $\euscS$ which at the same time is a locally free (resp.
free, coherent) $\cO_X$-module. A coherent right $\euscA$-module is a
sheaf of right $\euscA$-modules which at the same time is a coherent
$\cO_X$-module. 

\subsection{Coherent sheaves of right \boldmath{$\euscA$}-modules}
\label{sub:equiv-pathalgebra} 

Throughout \S \ref{sub:equiv-pathalgebra}, we assume that $Q$ is a
finite quiver, that is, $Q_0$ and $Q_1$ are both finite.  Let $M$ be
as in \S \ref{sub:quiverbundles}.

\bfsubsubsection{Twisted path algebra} \label{subsub:pathalgebra}

Let $\euscS=\oplus_{v\in Q_0}
\cO_X\cdot e_v$ be the free $\cO_X$-module generated by $Q_0$,
where $e_v$ are formal symbols, for $v\in Q_0$. We consider a
structure of commutative $\cO_X$-algebra on $\euscS$, defined by $e_v\cdot
e_{v'}=e_v$ if $v=v'$, and $e_v\cdot e_{v'}=0$ otherwise, for each
$v,v'\in Q_0$. Let  
$$
\euscM=\bigoplus_{a\in Q_1} M_a
$$
be a locally free sheaf of $\euscS$-bimodules, whose left (resp. right) 
$\euscS$-module structure is given by $e_v\cdot m=m$ if $m\in M_a$
and $v=ha$ (resp. $m\cdot e_v= m$ if $m\in M_a$ and $v=ta$), and
$e_v\cdot m= 0$ otherwise (resp. $m\cdot e_v= 0$ otherwise), for each
$v\in Q_0$, $a\in Q_1$, $m\in M_a$. The {\em $M$-twisted path algebra}
of $Q$ is the tensor $\euscS$-algebra of the $\euscS$-bimodule
$\euscM$, that is, 
$$ 
 \euscA=\bigoplus_{\ell\geq 0} \euscM^{\otimes_\euscS \ell}.
$$
Note that $\euscA$ is a locally free $\cO_X$-algebra. Furthermore, since
$Q$ is finite, $\euscA$ has a unit 
\begin{equation}\label{eq:unit}
1_\euscA=\oplus_{v\in Q_0} e_v.
\end{equation}

\bfsubsubsection{Coherent \boldmath{$\euscA$}-modules} 
\label{subsub:coh-A-mod}

We will show now that the category of $M$-twisted $Q$-sheaves is
equivalent to the category of coherent sheaves of right
$\euscA$-modules, or {\em coherent right $\euscA$-modules}. This
result is a direct generalisation of the corresponding equivalence of
categories for quiver modules 
(cf. e.g. \cite{ARS}). We define an equivalence functor from
the first to the second category. Let $\cR=(\cE,\phi)$ be an 
$M$-twisted $Q$-sheaf. Let $E=\oplus_{v\in Q_0} \cE_v$ as a coherent
$\cO_X$-module. The structure of right $\euscA$-module on $E$ is given
by a morphism of $\cO_X$-modules $\mu_\euscA: E\otimes_{\cO_X}\euscA\to E$
satisfying the usual axioms defining right modules over an algebra. 
Let $\pi_v:E\otimes_{\cO_X}\euscS=\oplus_{v,v'\in Q_0} \cE_v\otimes_{\cO_X}
\cO_X\cdot e_{v'} \to  \cE_v\otimes_{\cO_X}\cO_X\cdot e_v \cong
\cE_v$, be the canonical projection, and $\iota_v:\cE_v\hra E$ the
inclusion map, for each $v\in Q_0$. Let $\mu_v=\iota_v\circ\pi_v: E
\otimes_{\cO_X}\euscS \to E$. The morphism $\mu_\euscS=\sum_{v\in Q_0}
\mu_v :E\otimes_{\cO_X}\euscS \to E$ defines a structure of right 
$\euscS$-module on $E$. The tensor product of $E$ and $\euscM$ over
$\euscS$ is $E\otimes_\euscS \euscM\cong \otimes_{a\in Q_1}
E_{ta}\otimes_{\cO_X} M_a$; let $\pi_a: E\otimes_\euscS \euscM\to 
E_{ta}\otimes_{\cO_X} M_a$ be the canonical projection, for each $a\in
Q_1$. The morphism $\mu_{\euscM}=\sum_{a\in Q_1} \iota_{ha}\circ
\phi_a\circ\pi_a : E\otimes_\euscS \euscM \to E$ is a morphism of
$\euscS$-modules. Since $\euscA$ is the tensor $\euscS$-algebra of
$\euscM$, $\mu_{\euscM}$ induces a morphism of $\cO_X$-modules
$\mu_\euscA: E\otimes_{\cO_X}\euscA\to E$ defining a structure of right
$\euscA$-module on $E$. This defines the action of the equivalence
functor on the objects of the category of $M$-twisted $Q$-sheaves. 
It is straightforward to construct an action of the functor on
morphisms of $M$-twisted $Q$-sheaves, so this defines a functor
from the category of $M$-twisted $Q$-sheaves to the category of
coherent right $\euscA$-modules. We now define a functor from the
category of coherent right $\euscA$-modules to the category of
$M$-twisted $Q$-sheaves, and see that this new functor is an inverse
equivalence of the previous functor. Let $E$ be a coherent right
$\euscA$-module, with right $\euscA$-module structure morphism $\mu_\euscA:
E\otimes_{\cO_X}\euscA\to E$. The decomposition \eqref{eq:unit} is a
sum of orthogonal idempotents in 
$\euscA$ (i.e. $e^2_v=e_v$, $e_v\cdot e_{v'}=0$ for $v,v'\in Q_0$ with
$v\neq v'$), so $E=\oplus_{v\in Q_0}\cE_v$ with $\cE_v=\mu_\euscA
(E\otimes_{\cO_X} \cO_X\cdot e_v)\subset E$, for each $v\in Q_0$, and the
tensor product of $E$ and $\euscM$ over $\euscS$ is
$E\otimes_{\euscS}\euscM = \otimes_{a\in  Q_1} \cE_{ta}\otimes_{\cO_X}
M_a$. The restriction of 
$\mu_\euscA$ to $E\otimes_{\cO_X} \euscM$ induces a morphism of
$\euscS$-modules $\mu_\euscM:E\otimes_\euscS\euscM\to E$. The image of
$\cE_{ta}\otimes_{\cO_X} M_a$ under $\mu_\euscM$ is therefore in
$\cE_{ha}$, hence defines a morphism of $\cO_X$-modules
$\phi_a:\cE_{ta} \otimes _{\cO_X} M_a\to \cE_{ha}$, for each $a\in
Q_1$. This defines a functor from the category of coherent right 
$\euscA$-modules to the category of $M$-twisted $Q$-sheaves. It is
straightforward to define the action of this functor on 
morphisms and to prove that this functor, together with the
previous one, are inverse equivalences of categories (actually, of
$\cO_X$-categories, cf. e.g. \cite{ARS}). This completes the proof of
the following: 

\begin{proposition}\label{prop:equivalence-cat}
The category of coherent right $\euscA$-modules is equivalent to the
category of $M$-twisted $Q$-sheaves on $X$.
\end{proposition}

\subsection{Tensor products of stable twisted quiver bundles}
\label{sub:tensorprod}

As a simple application of Proposition \ref{prop:equivalence-cat}, we now
prove that the tensor product of two polystable  
twisted holomorphic quiver bundles is polystable as well. To do this, we
first define the appropriate notion of tensor product
of quiver sheaves. 
Let $Q=(Q_0,Q_1)$ and $Q'=(Q'_0,Q'_1)$ be two finite quivers
with the same vertex set $Q_0=Q_0'$, and tail and head maps
$t,h:Q_1\to Q_0$, $t',h':Q'_1\to Q'_0$, respectively. Let $M$
(resp. $M'$) be a collection of finite rank locally free sheaves
$M_a$ (resp. $M'_{a'}$) on $X$, for each $a\in Q_1$ (resp. $a'\in
Q_1'$). Let $\euscS=\oplus_{v\in Q_0} \cO_X\cdot e_v$ be a free sheaf
of $\cO_X$-algebras as in \S \ref{subsub:pathalgebra}.
Let $\euscM=\oplus_{a\in Q_1} M_a$, $\euscM'=\oplus_{a'\in Q'_1} M_{a'}$,
be locally free sheaves of $\euscS$-bimodules defined as in \S
\ref{subsub:pathalgebra}, and 
$$
\euscA =\bigoplus_{\ell=0}^{\infty} \euscM^{\otimes_\euscS \ell},\quad 
\euscA'=\bigoplus_{\ell=0}^{\infty} \euscM^{\prime\otimes_\euscS \ell},
$$
the $M$-twisted and $M'$-twisted path algebras of 
$Q$ and $Q'$, resp. Thus, the category of coherent right $\euscA$-modules
(resp. $\euscA'$-modules) is equivalent to the category of $M$-twisted
$Q$-sheaves (resp. $M'$-twisted $Q'$-sheaves) on $X$. Let $Q''$ be the quiver which has
the same vertices as $Q$ and $Q'$, and has the arrows of $Q$ and
$Q'$, i.e. $Q''=(Q''_0,Q''_1)$ is the quiver, with tail and head
maps $t'',h'':Q''_1\to Q''_0$, defined by 
$$
Q''_0=Q_0=Q'_0,~Q''_1=Q_1\coproduct Q_1',
$$
$$
t''a=ta, ~ h''a=ha \textnormal{~if~} a\in Q_1,
\textnormal{~and~} 
t''a'=t'a', ~ h''a'=h'a' \textnormal{~if~} a'\in Q'_1. 
$$ 
Let $M''$ be the collection of finite rank locally free sheaves
$M''_a$ on $X$, for each $a\in Q_1''$, given by $M''_{a}=M_{a}$ if
$a\in Q_1$ and $M''_{a'}=M'_{a'}$ if $a'\in Q_1'$. Let  
$$
\euscM''=\euscM\oplus\euscM'=\bigoplus_{a\in Q_0} M''_a.
$$
The $M''$-twisted path algebra of $Q''$ is 
$$
\euscA''=\bigoplus_{\ell=0}^{\infty}\euscM^{\prime\prime\otimes_\euscS \ell}\cong
\euscA\otimes_\euscS \euscA'.
$$
The category of coherent $\euscA''$-modules is equivalent to the
category of $M''$-twisted $Q''$-sheaves on $X$. Let now $E$
(resp. $E'$) be a coherent right $\euscA$-module (resp.
$\euscA'$-module). Since $\euscS$ is a commutative $\cO_X$-sheaf and
$E,E'$ are coherent $\euscS$-modules, their tensor product
$E''=E\otimes_\euscS E'$ is well defined and is again a coherent 
$\euscS$-module. We define the structure of a coherent right
$\euscA''$-module on 
$E''$ by the isomorphism $\euscA''\cong \euscA\otimes_\euscS \euscA'$: 
the action of $a\otimes a'\in (\euscA\otimes_\euscS \euscA')_x$ on
$e\otimes e'\in E''_x$, for each $x\in X$, is $(e\otimes e')\cdot(a\otimes a')
=e\cdot a\otimes e'\cdot a'$. 
Let now $\cR=(\cE,\phi)$ be the $M$-twisted
$Q$-sheaf corresponding to $E$, and $\cR'=(\cE',\phi')$ the 
$M'$-twisted $Q'$-sheaf corresponding to $E'$, by the equivalences
of categories of Propositin \ref{prop:equivalence-cat}. The
$M''$-twisted $Q''$-sheaf corresponding to their 
tensor product $E''$ is then $\cR''=(\cE'',\phi'')$, where
$\cE''_v=\cE_v\otimes_{\cO_X}\cE'_v$ for each $v\in Q_0$, and
$\phi''_{a}=\phi_{a}\otimes_{\cO_X}\id$ if $a\in Q_1$,
$\phi''_{a'}=\id\otimes_{\cO_X}\phi'_{a'}$ if $a'\in Q'_1$. Thus, $\cR''$ is
the {\em tensor product} of $\cR$ and $\cR'$.  

\begin{proposition} 
Let $\sigma,\tau,\tau'$ be collections of real numbers $\sigma_v$,
 $\tau_v$ and $\tau'_v$, respectively, with $\sigma_v>0$, for each
 $v\in Q_0$, and let $\tau''=\tau+\tau'$.
If $\cR$ is a $(\sigma,\tau)$-polystable holomorphic $M$-twisted
$Q$-bundle and $\cR'$ is a $(\sigma,\tau')$-polystable $M'$-twisted
holomorphic $Q'$-bundle, then their tensor product $\cR''$ is a
$(\sigma,\tau'')$-polystable holomorphic $M$-twisted $Q''$-bundle. 
\end{proposition}

\proof 
To define the vortex equations on $\cR$ and $\cR'$, resp., we fix a
family $q$ of hermitian metrics $q_a$ on $M_a$, for each $a\in Q_1$,
and a family $q'$ of hermitian metrics $q'_{a'}$ on $M'_{a'}$, for
each $a'\in Q_1'$, resp. 
By Theorem \ref{thm:HK-Q-bundles}, there is a hermitian metric
$H$ on $\cR$ satisfying the $(\sigma,\tau)$-vortex equations, and
a hermitian metric $H'$ on $\cR'$ satisfying the
$(\sigma,\tau')$-vortex equations. The Chern connection associated to
the metric $H''_v=H_v\otimes H'_v$ on $\cE''_v=\cE_v\otimes\cE'_v$,
for each $v\in Q_0$ with nonzero $\cE_v$ and $\cE'_v$, has curvature
$F_{H''_v}=F_{H_v}\otimes\id+\id\otimes F_{H'_v}$. It is now 
straightforward to prove that the collection $H''$ of hermitian
metrics $H''_v$ on $\cE''_v$, for each $v\in Q_0$, is a hermitian
metric on $\cR''$ satisfying the $(\sigma,\tau'')$-vortex equations.
Thus, $\cR''$ is\ $(\sigma,\tau'')$-polystable, again by Theorem
\ref{thm:HK-Q-bundles}. 
\qed

\section{Examples}
\label{sec:examples}

\subsection{Higgs bundles}
Let  $X$ be  a Riemann surface.  A \emph{Higgs bundle} on
$X$ is a pair $(E,\Phi)$, where $E$   is a holomorphic vector
bundle over $X$ and $\Phi \in H^0(\End(E) \otimes K)$ is a holomorphic
endomorphism of $E$ twisted by the canonical bundle $K$ of $X$.
The quiver here consists of one vertex and one arrow whose head and 
tail coincide and the twisting bundle is dual of the canonical line
bundle of $X$, i.e. the holomorphic tangent bundle $T'X$ of $X$. This
quiver, and the twisting bundle attached to its arrow, is represented
in Fig. 1. 

The Higgs bundle $(E,\Phi)$ is \emph{stable} if the usual slope stability 
condition
$\mu(E') < \mu(E)$ is satisfied for all proper $\Phi$-invariant
subbundles $E'$ of $E$.  The existence theorem of Hitchin and Simpson 
\cite{H,S} says that $(E,\Phi)$ is polystable if and only
if there exists  a hermitian metric $H$ on $E$ satisfying
\begin{equation}
  \label{eq:hitchin}
  \begin{aligned}
  F_H + [\Phi,\Phi^*] &= -\sqrt{-1}\mu \id_E \omega, 
  \end{aligned}
\end{equation}
where $\omega$ is the K\"ahler form on $X$, $\id_E$ is the
identity on $E$, and $\mu$ is a constant.  Note that taking the trace
in the first equation and integrating over $X$ we get $\mu = \mu(E)$.

There are many reasons why Higgs bundles are of interest, one of the
most important of which is the fact that there is a bijective
correspondence between isomorphism classes of poly-stable Higgs
bundles of degree zero on $X$ and isomorphism classes of semisimple
complex representations of the fundamental group of $X$. This
important fact is derived from a combination of the theorem of
Hitchin and Simpson mentioned above and an existence theorem for
equivariant harmonic metrics  proved by Donaldson \cite{D3}  and
Corlette \cite{C}. This correspondence can also be used to
study representations of $\pi_1 (X)$ in non-compact real Lie groups.  
In particular, by considering the group  $\U(p,q)$ one obtains
another interesting example of a twisted quiver bundle. To identify this 
quiver we observe that there
is a homeomorphism between the moduli space
of semisimple representation of $\pi_1(X)$ in ${\U(p,q)}$ and the moduli 
space of poly-stable zero degree Higgs bundles $(E,\Phi)$  of the form
\begin{equation}
  \label{upq-higgs-bundle}
  \begin{aligned}
  E &= V \oplus W, \\
  \Phi &=
  \left(
  \begin{smallmatrix}
    0 & \beta \\
    \gamma & 0
  \end{smallmatrix}
  \right), 
  \end{aligned}
\end{equation}
where $V$ and $W$ are holomorphic vector bundles on $X$ of rank $p$
and $q$, respectively, 
\begin{displaymath}
    \beta  \in H^0(\Hom(W,V) \otimes K)\quad \text{and} \quad 
    \gamma\in H^0(\Hom(V,W) \otimes K).
\end{displaymath}

The corresponding quiver, with the twisting bundle attached to each
arrow, is represented in Fig. 2. 
Now, for this twisted quiver bundle one can consider the general
quiver equations. Although they only  coincide with Hitchin's equations 
(\ref{eq:hitchin}) for a  particular choice  of the parameters, it turns
out that the other values are very important to study the topology
of the moduli of representations of $\pi_1(X)$ into $\U(p,q)$
\cite{BGG1}.

$$
\begin{array}{ccc}

\setlength{\unitlength}{.8mm}
\allinethickness{.5mm}
\begin{picture}(62,42)(0,0)
\put(31,21){\arc{24}{3.38}{9.2}}
\put(19,21){\circle*{2}}
\put(51,21){\makebox(0,0){\small{$T'X$}}}
\put(15,21){\makebox(0,0){$\mathcal{E}$}}
\put(39,36){\makebox(0,0){$\Phi$}}
\put(42.8,21){\vector(0,-1){1}}
\put(43,21){\vector(0,-1){1}}
\put(43.2,21){\vector(0,-1){1}}
\end{picture}

&
\qquad \qquad 
&

\setlength{\unitlength}{.8mm}
\allinethickness{.5mm}
\begin{picture}(62,42)(0,0)

\put(31,21){\arc{24}{9.65}{12.3}}
\put(31,21){\arc{24}{6.62}{9.2}}

\put(19,21){\circle*{2}} 
\put(43,21){\circle*{2}}

\put(21,7){\makebox(0,0){\small{$T'X$}}}
\put(21,35){\makebox(0,0){\small{$T'X$}}}

\put(15,21){\makebox(0,0){$V$}}
\put(48,21){\makebox(0,0){$W$}}

\put(39,6){\makebox(0,0){$\beta$}}
\put(39,36){\makebox(0,0){$\gamma$}}

\put(31,32.8){\vector(1,0){1}}
\put(31,33){\vector(1,0){1}}
\put(31,33.2){\vector(1,0){1}}

\put(31,8.8){\vector(-1,0){1}}
\put(31,9){\vector(-1,0){1}}
\put(31,9.2){\vector(-1,0){1}}

\end{picture}
\\ & & \\ 
\textnormal{Fig. 1.}
& & 
\textnormal{Fig. 2.}
\end{array}
$$

A very  important  tool to study topological properties of
Higgs bundle moduli spaces and hence moduli spaces of
representations of the fundamental group is to  consider the $\CC^*$-action 
on the moduli
space given by multiplying the Higgs field $\Phi$ by a non-zero
scalar.  A point $(E,\Phi)$ is a fixed
point of the $\CC^*$-action  action if and only if it is a variation of Hodge
structure, that is, 
\begin{equation}
  E = F_1 \oplus \cdots \oplus F_m
  \label{eq:variation-of-hodge}
\end{equation}
for holomorphic
vector bundles $F_i$ such that the restriction 
$$
\Phi_i := \Phi_{|F_i} \in H^0(\Hom(F_i,F_{i+1})\otimes K).
$$
A variation of Hodge structure is therefore a twisted quiver bundle,
whose twisting bundles are $M_a=T'X$, and the infinite quiver
represented in Fig. 3.
\begin{center}
\setlength{\unitlength}{.8mm}
\allinethickness{.5mm}
\begin{picture}(62,12)(0,0)
\multiput(0,0)(15,0){5}{\circle*{2}}
\multiput(-14,0)(15,0){6}{\line(1,0){13}}
\multiput(-5,0)(15,0){6}{{\vector(1,0){0}}}
\multiput(-11,1)(15,0){6}{{$\stackrel{T'X}{}$}}
\end{picture}
\end{center}\begin{center}
Fig. 3: Variations of Hodge structure. 
\end{center}

One can generalize the notion of Higgs bundle to consider
twistings by a line bundle other than the canonical bundle.
These have also very interesting geometry \cite{GR}.

\subsection{Quiver bundles and dimensional reduction}

Quiver bundles and their vortex  equations appear naturally in the 
context of dimensional reduction. To explain this, consider 
the manifold $X\times
G/P$, where $X$ is a compact \kahler\ manifold, $G$ is a connected
simply connected semisimple complex Lie group and $P\subset G$ is a
parabolic subgroup, i.e. $G/P$ is a flag manifold. The group $G$ (and
hence, its maximal compact subgroup $K\subset G$) 
act trivially  on $X$ and in  the standard way on $G/P$. The \kahler\
structure on $X$ together with a $K$-invariant \kahler\  structure on 
$G/P$ define a product \kahler\ structure on $X\times G/P$. 
We  can now consider a $G$-equivariant vector bundle over $X\times G/P$ and
study $K$-invariant solutions to the Hermitian--Einstein equations.
It turns that these invariant solutions correspond to special solutions
to the quiver vortex equations on a certain quiver bundle over $X$, where the
quiver is determined by the parabolic subgroup $P$.
In \cite{AG1} we studied the case in which $G/P=\PP^1$, the complex
projective line, which is obtained as the quotient of 
$G=\SL(2,\CC)$ by the subgroup of lower
triangular matrices, generalizing previous work by \cite{G1,G2,BG}.  
The general case has been  studied
in \cite{AG2}.
We will just mention here some of the main results  and refer the
reader to the above mentioned papers. 

A key fact is the existence of a quiver $Q$ with relations $\cK$
naturally associated to the subgroup $P$. A {\em relation} of the
quiver is a formal complex linear combination $r=\sum_j c_j p_j$ of
paths $p_j$ of the quiver (i.e. $c_j\in\CC)$, and a path in $Q$
is a sequence $p=a_0\cdots a_m$ of arrows $a_j\in Q_j$ which compose,
i.e. with $ta_{j-1}=ha_j$ for $1\leq j\leq m$:  
\begin{equation}\label{eq:path}
  p: \quad 
  {\bullet }\stackrel{a_m}{\lra }
  \bullet\stackrel{a_{m-1}}{\lra }\cdots\stackrel{a_0}{\lra}
  {\bullet } 
\end{equation}

The set of vertices of the quiver associated to $P$ coincides with the
set of irreducible representations of $P$. The arrows and relations
are obtained by studying certain isotopical decompositions related to
the nilradical of the Lie algebra of $P$. For example, for $\PP^1$,
$\PP^1\times \PP^1$ and $\PP^2$, the quiver is the disjoint union of
two copies of the quivers in Fig. 4, 5 and 6, respectively.
\begin{center}
\setlength{\unitlength}{.8mm}
\allinethickness{.5mm}
\begin{picture}(62,12)(0,0)
\multiput(0,0)(15,0){5}{\circle*{2}}
\multiput(-14,0)(15,0){6}{\line(1,0){13}}
\multiput(-5,0)(15,0){6}{{\vector(1,0){0}}}
\end{picture}
\end{center}\begin{center}
Fig. 4: $G/P=\PP^1$. 
\end{center}
$$
\begin{array}{ccc}
\setlength{\unitlength}{.8mm}
\begin{picture}(62,42)(-1,-1)
\thicklines
\matrixput(0,0)(10,0){7}(0,10){5}{\circle*{2}}
\matrixput(1,0)(10,0){6}(0,10){5}{\line(1,0){8}}
\matrixput(7,0)(10,0){6}(0,10){5}{{\vector(1,0){0}}}
\matrixput(4,-4)(10,0){1}(0,10){1}{{\scriptsize $a^{(1)}$}}
\matrixput(0,1)(10,0){7}(0,10){4}{\line(0,1){8}}
\matrixput(0,7)(10,0){7}(0,10){4}{{\vector(0,1){0}}}
\matrixput(-6,4)(10,0){1}(0,10){1}{{\scriptsize $a^{(2)}$}}
\matrixput(-1,0)(10,0){1}(0,10){5}{\line(-1,0){4}}
\matrixput(61,0)(10,0){1}(0,10){5}{\line(1,0){4}}
\matrixput(0,41)(10,0){7}(0,10){1}{\line(0,4){4}}
\matrixput(0,-1)(10,0){7}(0,10){1}{\line(0,-4){4}}
\end{picture}
&
\qquad \qquad 
&
\setlength{\unitlength}{0.8mm}
\begin{picture}(62,60)(-1,-1)
\thicklines
\multiput(0,0)(10,0){7}{\circle*{2}}
\multiput(10,10)(10,0){6}{\circle*{2}}
\multiput(20,20)(10,0){5}{\circle*{2}}
\multiput(30,30)(10,0){4}{\circle*{2}}
\multiput(40,40)(10,0){3}{\circle*{2}}
\multiput(50,50)(10,0){2}{\circle*{2}}
\multiput(1,0)(10,0){6}{\line(1,0){8}}
\multiput(11,10)(10,0){5}{\line(1,0){8}}
\multiput(21,20)(10,0){4}{\line(1,0){8}}
\multiput(31,30)(10,0){3}{\line(1,0){8}}
\multiput(41,40)(10,0){2}{\line(1,0){8}}
\multiput(51,50)(10,0){1}{\line(1,0){8}}
\multiput(7,0)(10,0){6}{{\vector(1,0){0}}}
\multiput(17,10)(10,0){5}{{\vector(1,0){0}}}
\multiput(27,20)(10,0){4}{{\vector(1,0){0}}}
\multiput(37,30)(10,0){3}{{\vector(1,0){0}}}
\multiput(47,40)(10,0){2}{{\vector(1,0){0}}}
\multiput(57,50)(10,0){1}{{\vector(1,0){0}}}
\multiput(14,-3.9)(10,0){1}{{\scriptsize $a^{(2)}$}}
\matrixput(10,1)(10,0){1}(0,10){1}{\line(0,1){8}}
\matrixput(20,1)(10,0){1}(0,10){2}{\line(0,1){8}}
\matrixput(30,1)(10,0){1}(0,10){3}{\line(0,1){8}}
\matrixput(40,1)(10,0){1}(0,10){4}{\line(0,1){8}}
\matrixput(50,1)(10,0){1}(0,10){5}{\line(0,1){8}}
\matrixput(60,1)(10,0){1}(0,10){5}{\line(0,1){8}}
\matrixput(10,7)(10,0){1}(0,10){1}{{\vector(0,1){0}}}
\matrixput(20,7)(10,0){1}(0,10){2}{{\vector(0,1){0}}}
\matrixput(30,7)(10,0){1}(0,10){3}{{\vector(0,1){0}}}
\matrixput(40,7)(10,0){1}(0,10){4}{{\vector(0,1){0}}}
\matrixput(50,7)(10,0){1}(0,10){5}{{\vector(0,1){0}}}
\matrixput(60,7)(10,0){1}(0,10){5}{{\vector(0,1){0}}}
\multiput(11.6,4)(10,0){1}{{\scriptsize $a^{(1)}$}}
\dottedline{2}(-5,-5)(55,55)
\matrixput(61,0)(10,0){1}(0,10){6}{\line(1,0){4}}
\matrixput(0,-1)(10,0){7}(0,10){1}{\line(0,-1){4}}
\matrixput(60,51)(10,0){1}(0,10){1}{\line(0,1){4}}
\end{picture}
\\ & & \\ 
\textnormal{Fig. 5:~} G/P=\PP^1\times\PP^1.
& & 
\textnormal{Fig. 6:~} G/P=\PP^2.
\end{array}
$$

In the case of the quiver associated to $\PP^1$, the set of relations
is empty, while for the quivers associated to $\PP^1\times\PP^1$ and
$\PP^2$, the relations $r_\lambda$ are given
by 
$$
r_\lambda=a^{(2)}_{\lambda- L_1}
a^{(1)}_\lambda-a^{(1)}_{\lambda- L_2} a^{(2)}_\lambda ,
$$
where $\lambda\in\ZZ^2$ is a vertex, $L_1$ and $L_2$ are the canoncial
basis of $\CC^2$, and $a^{(j)}_\lambda:\lambda\to \lambda-
L_j$ are the arrows going out from $\lambda$, for $j=1,2$. Given a set
$\cK$ of relations of the quiver $Q$, a holomorphic $(Q,\cK)$-bundle
(with no twisting bundles $M_a$) is defined as a holomorphic $Q$-bundle
$\cR=(\cE,\phi)$ which satisfies the relations $r=\sum_j c_j p_j$ in
$\cK$, i.e. such that $\sum_j c_j \phi(p_j)=0$, where
$\phi(p):\cE_{ta_m}\to\cE_{ha_0}$ is defined for any path \eqref{eq:path}
as the composition $\phi(p):=\phi_{a_0}\circ\cdots\circ\phi_{a_m}$. 

Let $(Q,\cK)$ be the quiver with relations associated to $P$. One has
an equivalence of categories  
$$
\left\{\begin{array}{c}
{\rm coherent~} G {\rm -equivariant} \\
{\rm sheaves~on~} X\times G/P
\end{array}\right\}
\longleftrightarrow
\left\{\begin{array}{c}
(Q,\cK){\rm -sheaves~on~} X
\end{array}\right\}.
$$
The holomorphic $G$-equivariant vector bundles on $X\times G/P$ and the
holomorphic $(Q,\cK)$-bundles on $X$ are in correspondence by this
equivalence. Thus, the category of $G$-equivariant holomorphic vector
bundles on $X\times (\PP^1)^2$ and $X\times \PP^2$ is equivalent to
the category of commutative diagrams of holomorphic quiver bundles on
$X$ for the corresponding quiver $Q$. 
 If we now fix a total order in the set of vertices, 
any coherent $G$-equivariant sheaf 
$\cF$ on $X\times G/P$ admits a $G$-equivariant sheaf filtration
\begin{equation}\begin{gathered}\label{eq:equi-hol-fil}
\bcF:\, \holfil ,\\
\cF_s/\cF_{s-1}\cong p^*\cE_{\lambda_s}\otimes q^*\cO_{\lambda_s},
\quad 0\leq s\leq m,
\end{gathered}\end{equation}
where $\{\lambda_0,\lambda_1,\ldots ,\lambda_m\}$ is a
finite subset of vertices, listed in ascending order,
 $\cE_0,\ldots,\cE_m$ are
non-zero coherent sheaves on $X$  with trivial $G$-action,
and $\cO_{\lambda_s}$ is the homogeneous bundle over $G/P$ corresponding
to the representation $\lambda_s$. 
The maps $p$ and $q$ are the canonical projections from
$X\times G/P$ to $X$ and $G/P$, respectively.
If $\cF$ is a holomorphic $G$-equivariant vector bundle, 
then $\cE_0,\ldots,\cE_m$ are holomorphic vector bundles.

The  appropriate
equation to consider on a filtered bundle \cite{AG1}  is a
 deformation of  Hermite--Einstein equation which  involves as
many parameters $\tau_0,\tau_1,\ldots ,\tau_m\in\RR$ as steps are in
the filtration, and has the form  
\begin{equation}\label{eq:tau-metric-HEE}
\imag\Lambda F_h=\begin{pmatrix}
                \tau_0 I_0      &               &        &              \\ 
                                & \tau_1 I_1    &        &              \\ 
                                &               & \ddots &              \\ 
                                &               &        & \tau_m I_m
\end{pmatrix},
\end{equation}
where the RHS is a diagonal matrix, written in blocks corresponding
to the splitting which a hermitian metric $h$ defines in the
filtration $\bcF$. If $\tau_0=\cdots=\tau_m$, then
\eqref{eq:tau-metric-HEE} reduces to the Hermite--Einstein equation. 
As in the ordinary Hermite--Einstein equation, the existence of
invariant solutions to the $\tau$-Hermite--Einstein equation on an
equivariant holomorphic filtration is related to a stability condition
for the equivariant holomorphic filtration
which naturally involves the parameters.

Let $\cF$ be a $G$-equivariant holomorphic vector bundle on $X\times
G/P$. Let  $\bcF$ be the $G$-equivariant holomorphic filtration associated to  $\cF$
 and
$\cR=(\cE ,\phi)$ be its corresponding  holomorphic $(Q,\cK)$-bundle on
$X$, where $(Q,\cK)$ is  the quiver with relations associated to $P$.
Then $\bcF$ has a $K$-invariant solution to the  $\tau$-deformed 
Hermite--Einstein equations 
if and only if the vector bundles $\cE_\lambda$ in $\cR$
admit hermitian metrics $H_\lambda$ on $\cE_\lambda$,
for each vertex  $\lambda$ with $\cE_\lambda\neq 0$, satisfying
\begin{equation}\label{eq:vortex-connection}
\imag n_\lambda \Lambda F_{H_\lambda}
+\sum_{a\in h^{-1}(\lambda)}\phi_a\circ\phi_a^*-\sum_{a\in t^{-1}(\lambda)}\phi_a^*\circ\phi_a 
=\tau_\lambda'\id_{\cE_\lambda}, \quad\quad 
\end{equation}
where $n_\lambda$  is the multiplicity of the irreducible
representation corresponding to the vertex  $\lambda$
and $\tau_\lambda'$ are related to $\tau_\lambda$ by the choice of the
$K$-invariant metric on $G/P$.
It is not difficult to show that the stability of the filtration 
coincides with the stability of the quiver bundle
where the parameters $\sigma_\lambda$ in the general stability condition
for a quiver bundle equal the integers $n_\lambda$. This, together
with the dimensional reduction obtainment of the equations, provides
with an alternative  proof of the Hitchin--Kobayashi correspondence
for these special quiver bundles.

Although the quiver bundles obtained by dimensional reduction
on $X\times G/P$ are not
twisted, it seems  that twisting may appear if one considers
 dimensional reduction on more general $G$-manifolds ---
this is something to  which we plan to come back in the future.


\bfsubsection*{Acknowledgements} 

This research has been partially supported by
the Spanish MEC under the grant PB98--0112.
The research of L.A. was partially supported by the Comunidad
Aut\'onoma de Madrid (Spain) under a FPI Grant, and by a UE Marie Curie
Fellowship (MCFI-2001-00308).
The authors are members of VBAC (Vector Bundles on Algebraic Curves),
which is partially supported by EAGER (EC FP5 Contract no. HPRN-CT-2000-00099)
and by EDGE (EC FP5 Contract no. HPRN-CT-2000-00101).
We also want to thank the Erwin Schr\"odinger International Institute 
for Mathematical Physics for the hospitality and the support during 
the final preparation of the paper. 


\end{document}